\documentclass[11pt, psamsfonts]{amsart}
\usepackage{amssymb, amsmath, amsthm}
\usepackage{graphicx}
\usepackage[final, hypertex]{hyperref}
\usepackage[a4paper, centering]{geometry}
\newcommand{\picturename}[1]{#1.eps}

\newtheorem{theorem}{Theorem}[section]
\newtheorem{prop}[theorem]{Proposition}
\newtheorem{lemma}[theorem]{Lemma}

\theoremstyle{definition}
\newtheorem{dhef}[theorem]{Definition}
\theoremstyle{remark}
\newtheorem{rk}[theorem]{Remark}

\newtheorem{claim}[theorem]{{\it Claim}}
\newenvironment{claim-b}
   {\begin{claim}\rm}
   {\end{claim}}

\long\def\elimina#1{} 

\def\R{\mathbb{R}}

\def\N{\mathbb{N}}

\def\pscal#1#2{\left\langle#1,\,#2\right\rangle}

\def\haus{\mathcal{H}^{n-1}}

\def\gau#1{\rho_{#1}}
\def\gauge{\rho}
\def\gaugem{\rho_{-}}
\def\pgauge{\gauge^0}

\def\distb#1{d_{#1}}
\def\dist{\distb{\Omega}}
\def\distm{\distb{\Omega}^-}

\def\curvg{\tilde{\kappa}}
\def\curvgm{\tilde{\kappa}^-}
\def\nor{\nu}

\def\len{l}

\def\Lip{\textrm{Lip}}
\def\Lipr{\Lip_{\rho}}
\def\bw{\overline{W}}

\def\vf{v}

\def\uz{u}
\def\Cb{C_b}
\def\dato{\overline{u}}
\def\dt{\delta t^n_i}
\def\macroy{y}
\def\nt{\vec{n}}
\def\normp#1{{\left\|#1\right\|}_p}
\def\Uz{\vec{U}_0}
\def\Usol{\vec{U}}
\def\cz{\vec{c}_0}
\def\gd{\gauge_{\Delta}}
\def\elle#1{L^{#1}(A,\R^3)}
\def\wp#1{W^{1,#1}(A,\R^3)}

\DeclareMathOperator{\inte}{int} \DeclareMathOperator{\spt}{supp}
 \DeclareMathOperator{\dive}{div}
 
\DeclareMathOperator{\proj}{\Pi} \DeclareMathOperator{\curl}{curl}

\begin{document}
\title[Electrodynamics of hard superconductors]
{A variational approach to the macroscopic electrodynamics
of anisotropic hard superconductors}%

\author[G.~Crasta]{Graziano Crasta}
\address{Dipartimento di Matematica ``G.\ Castelnuovo'', Univ.\ di Roma I\\
P.le A.\ Moro 2 -- 00185 Roma (Italy)}
\email[Graziano Crasta]{crasta@mat.uniroma1.it}

\author[A.~Malusa]{Annalisa Malusa}
\email[Annalisa Malusa]{malusa@mat.uniroma1.it}

\date{December 20, 2006}

\keywords{Minimum problems with constraints, Euler equation, hard superconductors, Bean's model}
\subjclass[2000]{Primary 35C15; Secondary 49J30, 49J45, 49K20}

\begin{abstract}
We consider the Bean's critical state model for anisotropic superconductors.
A variational problem solved by the quasi--static evolution
of the internal magnetic field is obtained as the $\Gamma$-limit of functionals
arising from the Maxwell's equations combined with a power law for the
dissipation. Moreover, the quasi--static approximation of the internal
electric field is recovered, using a first order necessary condition.

If the sample is a long cylinder subjected to an axial uniform external
field, the macroscopic electrodynamics is explicitly determined.
\end{abstract}

\maketitle

\section{Introduction}


It is well known that a superconductor is a conductor which is able to pass an electric current
without dissipation. The transition from the normally conducting state
to the superconducting one occurs at a critical temperature $T_c$, depending on the material,
below which the material exhibits (almost) perfect conductivity.
We are interested in the response of a superconducting material  to an applied external
magnetic field $\vec{H}_s$  under isothermal conditions below its critical temperature.

The superconductors can be classified in terms of a material parameter $\kappa>0$,
known as the Ginzburg--Landau parameter.
For the so called type--II superconductors (corresponding to $\kappa>1/\sqrt{2}$)
there exist two critical magnetic field intensities,
${H}_{c_1} < {H}_{c_2}$, such
that for $|\vec{H}_s| < {H}_{c_1}$ the material is in the superconducting state and the
magnetic field is excluded from the bulk of the sample except in thin boundary layers,
while for $|\vec{H}_s| > {H}_{c_2}$ the material
behaves as a normal conductor,
and the magnetic field penetrates it fully.
For ${H}_{c_1}<|\vec{H}_s|< {H}_{c_2}$ a third
state exists, known as ``mixed state'' as well as ``vortex state''.
The mixed state is characterized by a partial penetration
of the magnetic field into the sample, which occurs, at a mesoscopic level, by
means of thin filaments of normally conducting material carrying magnetic flux
and circled by a vortex of superconducting current.
Type--I superconductors are those with ${H}_{c_1}={H}_{c_2}$, so that
the mixed state does not occur.

Although from a physical point of view the most interesting description of the mixed
state for type--II superconductors is given by the mesoscopic Ginzburg--Landau model,
in designing magnets and other large--scale applications of superconducting materials,
engineers use macroscopic models, involving averaged variables. One of the
most reliable macroscopic models is the Bean's critical state model (see \cite{Bean}).
We refer to \cite{Cha} for a derivation of this model
as a macroscopic version of the Ginzburg--Landau model
under suitable assumptions.

The basic idea of the Bean's phenomenological model for isotropic materials is that,
because of physical limitations imposed by the material properties, the current
density $|\vec{J}|$ cannot
exceed a critical value $J_c$ without destroying the superconducting phase.

Moreover it is assumed that any electromotive
force due to external field variations induces the maximum current density flow,
according to the most effective way of shielding field variations. Hence
the current density $|\vec{J}|$ is forced to be $J_c$ in the part of the sample
where the field is penetrated, while $\vec{J}=0$ in the remaining part of the sample.

At a mesoscopic level, the critical current density $J_c$ corresponds to the
balance between a repulsive vortex--vortex interaction and attractive forces
towards the pinning centers. At a macroscopic level, the electric field is zero
when $|\vec{J}|<J_c$ and abruptly rises to arbitrarily large values if $J_c$ is
overrun.

Since for isotropic materials it is well--established that the electric field and the
current have the same direction (at least for slowly varying external fields),
in the original Bean's model
the Ohm's law is replaced by a vertical current--voltage law
\[
\vec{E} \parallel \vec{J}\,, \quad
|\vec{J}|\leq J_c\,,\qquad
\vec{E}= 0 \ \textrm{if}\ |\vec{J}| < J_c\,,
\]
which can be interpreted as the limit
for $p\to +\infty$ in the power--law
\begin{equation}\label{f:power}
\vec{E} \parallel \vec{J}\,, \quad |\vec{E}|=e_c \left(\frac{|\vec{J}|}{|\vec{J}_c|}\right)^p
\end{equation}
(see, e.g., \cite{Bran,BaPr,BLc}).
Notice that the direction of the electric field is obtained by exploiting the isotropy of the
material and cannot be obtained as a consequence of the power law approximation, which involves
only the intensities of the fields.
Hence this approach is not appropriate to deal with anisotropic materials
(see anyhow \cite{BLb} for some
explicit computation of the electric field in the case of infinite slab geometry,
and \cite{BKR} in the case of a cylindrical body with elliptic section).
On the other hand, many type--II superconducting materials are anisotropic. Moreover,
even if the material is isotropic, the presence of flux-flow Hall effects has to be described
in terms of an anisotropic resistivity (see, e.g., \cite{Bran}).

Concerning the constraints on the current, in the anisotropic case the Bean's law dictates that
there exists a compact set $\Delta$ containing the origin as an interior point, such that
$\vec{J}$ cannot lie outside $\Delta$ without destroying the superconducting state.
Moreover $\vec{J}\in \partial \Delta$ in the penetrated region, and $\vec{J}=0$ elsewhere.

Section \ref{s:phis} of this paper is devoted to the description of the Bean's law as
a limit of a power--like law fulfilled by the dissipation. This approach
allows us to determine \textit{a priori} the direction of the electric field in terms of
the direction of the current for anisotropic materials.

In Section \ref{s:bean} we deal with a variational statement of the critical state
proposed in \cite{BLd}, which takes
the form of a quasistatic evolution of the penetrated magnetic field, obtained combining
the finite--difference expression of Faraday's law with the Bean's law.
We shall give a mathematical justification of the variational model
as a limiting case of the power law model for dissipation.
The result is proposed in terms of
$\Gamma$--convergence of functionals, which is nowadays a classical tool in the mathematical methods for
the material science (see, e.g., \cite{Brai,DM} and the references therein).

In the last part of the paper we focus our attention to the special case of a long
cylindrical anisotropic superconductor placed into a nonstationary, uniform axial
magnetic field. The parallel geometry enables us to make a two dimensional reduction
of the problem,
which can then be explicitly solved.

The plan of this part of the paper is the following. In Section \ref{s:bean2} we describe the
two dimensional reduction of the problem in the case of parallel geometry.
In Section \ref{s:distance} we introduce some
technical tool needed for the analytical description of the fields.
In Section \ref{s:minp} we find explicitly the solution to a general class of minimum problems
with a gradient constraint and we determine the Euler equation solved by the optimal function
coupled with its dual function. Moreover we find the explicit form of the dual function.
In Section \ref{s:evol} the previous results are applied to the variational model for the critical
state, and we find the quasistatic evolution of both the magnetic field and the dissipation
inside the superconductor. Finally, a passage to the limit on the time layer gives
the explicit form of the macroscopic electrodynamics.

For what concerns the magnetic field, our result
generalizes the one, valid for isotropic materials, obtained by Barrett and Prigozhin
in \cite{BaPr} with a different method
based on a evolutionary variational inequality.
Due to the fact that the magnetic field is explicitly known,
we can compute the full penetration time in the case of monotonic external fields
as well as we can depict the well known hysteresis phenomenon.
On the other hand, the knowledge of the electric field allows us to give a detailed
description of the evolution of the dissipation
for cylinders of anisotropic materials with a general geometry of the cross section
(see Figures \ref{BB3} and \ref{BB4}).

\section{Notation and preliminaries}

For $\xi\in \R^N$, $|\xi|$ will be the Euclidean norm, and
$\pscal{\xi}{\xi'}$ will denote the scalar product with $\xi'\in \R^N$.
The symbol $\times$ will be used for the cross product of vectors
in $\R^3$.
Given $a,b\in\R$, $a\vee b$ and $a\wedge b$ will denote
respectively the maximum and the minimum of $a$ and $b$.

\smallskip

Given $A\subset\R^N$, we shall denote by $\Lip(A)$, $C(A)$,
$\Cb(A)$ and $C^k(A)$, $k\in\N$ the set of functions $u\colon
A\to\R$ that are respectively Lipschitz continuous, continuous,
bounded and continuous, and $k$-times continuously differentiable
in $A$. Moreover, $C^{\infty}(A)$ will denote the set of functions
of class $C^k(A)$ for every $k\in\N$, while $C^{k,\alpha}(A)$ will
be the set of functions of class $C^k(A)$ with H\"older continuous
$k$-th partial derivatives with exponent $\alpha\in [0,1]$.
Finally $L^p(A)$, $W^{1,p}(A)$, $W^{1,p}_0(A)$, and $L^p(A,\R^d)$,
$W^{1,p}(A,\R^d)$, $W^{1,p}_0(A,\R^d)$, $d>1$, will be the usual Lebesgue
and Sobolev spaces of scalar or vectorial functions respectively.

\smallskip

Let $X$ be subset of $\R^N$, $N\geq 2$. We shall denote by
$\partial X$ its boundary, by $\inte X$ its interior, and by $\overline{X}$ its closure.
The characteristic function of $X$
will be denoted by $\chi_X$.
The set $X$ is said to be
of class $C^k$, $k\in\N$,
if for every point $x_0\in\partial X$
there exists a ball $B=B_r(x_0)$ and a one-to-one
mapping $\psi\colon B\to D$ such that
$\psi\in C^k(B)$, $\psi^{-1}\in C^k(D)$,
$\psi(B\cap X)\subseteq\{x\in\R^n;\ x_n > 0\}$,
$\psi(B\cap\partial X)\subseteq\{x\in\R^n;\ x_n = 0\}$.
If the maps $\psi$ and $\psi^{-1}$ are of class
$C^{\infty}$ or $C^{k,\alpha}$ ($k\in\N$, $\alpha\in [0,1]$),
then $X$ is said to be of class
$C^{\infty}$ or $C^{k,\alpha}$ respectively.

\smallskip

Let $D\subset\R^N$ be a compact convex set containing
$0$ as an interior point and with boundary of class $C^2$.
The gauge function $\gau{D}$ of the set $D$ is the convex, positively
1--homogeneous function defined by
\[
\gau{D}(\xi) = \inf\{ t\geq 0;\ \xi\in t D\}\,,
\qquad\xi\in\R^N.
\]
Since $D$ is a compact set containing a neighborhood of
$0$, there exist two positive constants $c_1<c_2$ such that
\begin{equation}\label{robdd}
c_1|\xi| \leq \gau{D}(\xi)\leq c_2 |\xi|\,, \qquad \forall\ \xi \in \R^N\,.
\end{equation}

The indicator function of the set $D$ is defined
by
\begin{equation}\label{indk}
I_{D}(\xi)=
\begin{cases}
0 & \textrm{if}\ \xi\in D \\
+\infty & \textrm{if}\ \xi\not\in D\,.
\end{cases}
\end{equation}

Since $D$ has a smooth boundary, the subgradient of the indicator function
can be explicitly computed, and
\begin{equation}\label{subindk}
\partial I_{D}(\xi)=
\begin{cases}
\{\alpha D\gau{D}(\xi)\colon \alpha \geq 0\} & \textrm{if}\ \xi\in \partial D \\
\emptyset & \textrm{if}\ \xi\not\in D \\
\{0\} & \textrm{if}\ \xi\in \inte D
\end{cases}
\end{equation}
(see e.g.\ \cite{Rock}, Section 23).

In what follows, a family of objects, even if parameterized by a continuous
parameter, will also be often called a ``sequence'' not to overburden
notation.

\section{The physical setting}\label{s:phis}

Since we shall deal with a macroscopic model,
the coarse--grained electrodynamics will be formulated
in terms of
\begin{itemize}
\item[(i)] the flux density $\vec{B}$ within the sample, which is the average of the
microscopic field intensity;
\item[(ii)] the magnetic field $\vec{H}$, which is assumed to be linearly connected
to $\vec{B}$ by $\vec{B}=\mu_0 \vec{H}$;
\item[(iii)] the averaged current density $\vec{J}$, linked to  $\vec{H}$ by the
Amp\`ere's law $\vec{J}=\curl \vec{H}$.
\end{itemize}
Moreover, on neglecting finite size effects, we can assume that
\begin{itemize}
\item[(iv)] the magnetic source $\vec{H}_s$ enters
as a boundary condition for the flux density at the surface of the sample,
requiring that $\nt\times\vec{H} = \nt\times\vec{H}_s$ on that surface,
where $\nt$ is the outward normal vector.
\end{itemize}

Time variations of the external field $\vec{H}_s$ induce in the sample an electric field
$\vec{E}$, according to Faraday's
law.
On the other hand,
the electric field leads
a current $\vec{J}$ which
induces an internal magnetic field $\vec{H}$,
according to Amp\`ere's law.
Finally, we recall that $\dive \vec{H}=0$ (Gauss' law).

Summarizing, we are considering the so-called eddy current model
for Maxwell equations:
\[
\begin{cases}
\curl\vec{E}+\mu_0\dfrac{\partial \vec{H}}{\partial t}=0
&\textrm{(Faraday's law)},\\
\vec{J}=\curl \vec{H}
&\textrm{(Amp\`ere's law)},\\
\dive\vec{H} = 0
&\textrm{(Gauss' law)},\\
\nt\times(\vec{H} -\vec{H}_s) = 0
&\textrm{on the surface}.
\end{cases}
\]

In order to complete the physical setting,
it remains to find an appropriate version of the constitutive law
$\vec{E}(\vec{J})$. Our starting point is a reading of the macroscopic behaviour
of the material in terms of the dissipation
$\dot{S}=\pscal{\vec{E}(\vec{J})}{\vec{J}}$
(here $S$ denotes the entropy of the system).
Namely, the Bean's model dictates that if $\vec{J}$ is in the interior of
an allowed region $\Delta$, which is assumed to be a convex compact set of $\R^3$
having the origin as an interior point,
then $\dot{S}$ vanishes as no electric field is generated in stationary
condition, while when $\vec{J}$ touches
the boundary $\partial \Delta$ of the allowed region a huge dissipation occurs,
destroying the superconducting phase.
The fact that $\dot{S}=0$ if $\vec{J}$ is an interior point of $\Delta$,
while $\dot{S} = +\infty$ if $\vec{J}\not\in\Delta$
suggests that the Bean's law can be
interpreted as the limit as $p\to\infty$ of a power law for the dissipation
$\pscal{\vec{E}_p(\vec{J})}{\vec{J}}$, that is
\begin{equation}\label{f:pdissip}
\pscal{\vec{E}_p(\vec{J})}{\vec{J}}=\frac{c}{p}\left(\gauge_\Delta(\vec{J})\right)^{p}\,.
\end{equation}
This approximation allows us to recover the direction of
the electric field $\vec{E}(\vec{J})$ when
$\vec{J}\in\partial\Delta$.
Namely, we claim that, under the physically consistent assumption that
\begin{equation}\label{f:nulel}
\lim_{t\to 0^+}\vec{E}_p(t\vec{J})=0\,, \quad \forall\ \vec{J}\in\R^3\,,
\end{equation}
the relation (\ref{f:pdissip}) implies, for $p>1$,
\begin{equation}\label{f:pelet}
\vec{E}_p(\vec{J})=\frac{c}{p}\left(\gauge_\Delta(\vec{J})\right)^{p-1} D \gauge_\Delta(\vec{J})\,.
\end{equation}
In order to prove (\ref{f:pelet}), we
differentiate (\ref{f:pdissip}), obtaining
\[
\pscal{D\vec{E}_p(\vec{J})}{\vec{J}}+ \vec{E}_p(\vec{J})=
c\left(\gauge_\Delta(\vec{J})\right)^{p-1} D \gauge_\Delta(\vec{J})\,.
\]
Hence, fixed $\vec{J}\neq 0$, the function $v(t)= \vec{E}_p(t\vec{J})$, $t>0$,
is a solution of the O.D.E.
\begin{equation*}
\begin{cases}
t\, v'(t)+v(t)=c\, t^{p-1}w \,, \\
v(1)= \vec{E}_p(\vec{J})\,,
\end{cases}
\end{equation*}
where $w=\left(\gauge_\Delta(\vec{J})\right)^{p-1} D \gauge_\Delta(\vec{J})$. Then
\[
v(t)=\frac{1}{t}\left[\vec{E}_p(\vec{J})+\frac{cw}{p}(t^p-1)\right]\,,
\]
and, by (\ref{f:nulel}),
\[
\vec{E}_p(\vec{J})-\frac{c}{p}\left(\gauge_\Delta(\vec{J})\right)^{p-1} D \gauge_\Delta(\vec{J}) =
\lim_{t\to 0^+}[\vec{E}_p(\vec{J})+\frac{cw}{p}(t^p-1)]=0\,.
\]

As a consequence of (\ref{f:pelet}), the direction of $\vec{E}_p(\vec{J})$ is given
by $D \gauge_\Delta(\vec{J})$, and it does not depend on $p$.
In conclusion, recalling (\ref{subindk}),
the anisotropic version of the Bean's law is the following.
\begin{itemize}
\item[(B1)] There exists a convex
compact set $\Delta$ containing the origin as an interior point and
such that $\vec{J}(x,t)\in \Delta$ for every $x\in\Omega$, $t\geq 0$.
\item[(B2)] If we denote by $I_\Delta$ the indicator function
of the set $\Delta$, the constitutive law $\vec{E}(\vec{J})$ is given by
$\vec{E}\in \partial I_\Delta(\vec{J})$.
\end{itemize}


\section{A variational model for the mixed state}\label{s:bean}

Let $A\subset\R^3$ be the region occupied by the superconductor, and
$\partial A$ its surface.
We assume that $A$ is a bounded, open subset of $\R^3$, and that
$\partial A$ is of class $C^{1,1}$.
Moreover, we assume that $A$ has no enclosed cavity.
In our model the presence of cavities is not allowed,
since the boundary condition in such a cavity cannot be given in terms of
the external magnetic field.
(The magnetic field must be constant on the boundary of every enclosed cavity.)
In order to avoid the presence of cavities, we assume that $A$ has second Betti number $0$
(see e.g.\ \cite{Gold} for the definition of Betti numbers).

Let $T>0$ be fixed. For every $n\in\N^+$, we take a partition
$P^n=\left\{t^n_i\right\}_{i=0}^{k(n)}$
of the interval $[0,T]$, and we set $\delta t^n_i=t^n_i-t^n_{i-1}$.
The least action principle proposed in \cite{BLd} is the following: starting from an initial field profile
$\vec{H}^n_i(x)$ in $A$ at the time layer $i$,
and under a small change of the external drive, the new profile $\vec{H}^n_{i+1}(x)$
at the time layer $i+1$
is the unique solution to the minimum problem
\begin{equation}\label{varmod3d}
\min\left\{
\int_A  |\vec{V}-\vec{H}^n_i|^2\, dx;\
\vec{V}\in \vec{H}_s(t^n_{i+1})+X_2,\
\curl\vec{V}\in \Delta\ \textrm{a.e.}
\right\}\,,
\end{equation}
where, for every $p>1$,
\begin{equation*}
X_p = \left\{\vec{V}\in \elle{p};\
\curl\vec{V}\in \elle{p},\,
\dive\vec{V} = 0,\,
\nt\times\vec{V} = 0\, \textrm{on}\, \partial A
\right\}.
\end{equation*}

In \cite{BLd}
the trustworthiness of this variational principle was motivated
in terms of nonequilibrium thermodynamical principles
in analogy with the ohmic case.
In this section we shall give a mathematical justification, in terms of
$\Gamma$--convergence of functionals, of the variational model
as a limiting case of the power law model for dissipation.

Our starting point is the fact that, assuming that the power law (\ref{f:pdissip})
holds true, the discretized version of the Faraday's law
\begin{equation}\label{f:farlp}
\curl\vec{E}^n_{i+1}+\frac{\mu_0}{c\dt}(\vec{H}^n_{i+1}-\vec{H}^n_{i})=0
\end{equation}
is the Euler equation for the functional
\[
F_p(\vec{V})=\int_A \frac{1}{p}\gd(\curl\vec{V})^p+
\lambda\left|\vec{V}-\vec{H}^n_i\right|^2\, dx\,,\qquad
\vec{V}\in \vec{H}_s(t^n_{i+1})+X_p
\]
where $\lambda := \frac{\mu_0}{2c\dt}>0$.

In other words, it can be easily checked that the magnetic field $\vec{H}^n_{i+1}$
is the unique minimum point of $F_p$ in $\vec{H}_s(t^n_{i+1})+ X_p$.
It is clear that, up to a translation of a constant vector,
it is equivalent to minimize $F_p$ on $X_p$ by
changing $\vec{H}^n_i$ with $\vec{H}^n_i-\vec{H}_s(t^n_{i+1})$.

The following properties of the spaces $X_p$ will be useful in the sequel
(see \cite[Thm.~2.2]{YLZ} and \cite[Thm.~3.1]{Wa}).

\begin{theorem}\label{t:Xp}
Let $A\subset\R^3$ be a bounded open set of class $C^{1,1}$.
Let $1<p<\infty$. Then the following hold.
\begin{itemize}
\item[(i)]
Every function $\vec{V}\in X_p$ belongs to $\wp{p}$.
Moreover, if $p>3$, then $\vec{V}\in C^{0,\alpha}(\overline{A})$
and there exists a constant $C_p>0$, depending only on $p$ and $A$,
such that
\[
{\left\|\vec{V}\right\|}_{C^{0,\alpha}(\overline{A})}\leq
C_p \normp{\curl\vec{V}}\,.
\]

\item[(ii)]
If $A$ has second Betti number $0$, then there exists a constant
$M_p > 0$, depending only on $p$ and $A$,
such that
\[
\normp{D\vec{V}} \leq M_p \normp{\curl\vec{V}}
\]
for every $\vec{V}\in X_p$.
\end{itemize}
\end{theorem}

\begin{rk}\label{r:Xp}
{}From Theorem~\ref{t:Xp} we infer that $X_p$ is a (proper)
subspace of $\wp{p}$.
Moreover, for $p>3$ the map defined by
$X_p\ni\vec{V}\mapsto\normp{\curl\vec{V}}$ is a norm in $X_p$ equivalent
to the standard $W^{1,p}$ norm.
\end{rk}

From now on we shall always assume $p>3$. In this case,
by Theorem \ref{t:Xp}(i), a function $\vec{V}\in X_p$ is continuous
in $\overline{A}$,
and the boundary condition $\nt\times \vec{V}$ is understood to
be pointwise fulfilled. On the other hand, the divergence--free requirement
on $\vec{V}\in X_p$ is understood in the sense of distributions, i.e.
\[
\int_A \vec{V}\cdot D\vec{\Psi} =0\,, \qquad \forall \vec{\Psi}\in C^\infty_0(A, \R^3)\,.
\]

For our subsequent considerations, we need to define all the functionals $F_p$
on the same Banach space.
For this reason,
naming $\vec{U}_0=\vec{H}^n_i-\vec{H}_s(t^n_{i+1})$,
we set
\begin{equation}\label{f:Fp}
F_p(\vec{V})=
\begin{cases}
\displaystyle
\int_A \frac{1}{p}\gd(\curl\vec{V})^p+
\lambda\left|\vec{V}-\vec{U}_0\right|^2\, dx\,,
&\textrm{if $\vec{V}\in X_p$},\\
+\infty,
&\textrm{otherwise in $\elle{2}$}.
\end{cases}
\end{equation}
We also define the functional
\begin{equation}\label{f:F}
F(\vec{V})=
\begin{cases}
\displaystyle\int_A
I_{\Delta}(\curl\vec{V})+
|\vec{V}-\Uz|^2\, dx\,,
&\textrm{if $\vec{V}\in X_2$},\\
+\infty,
&\textrm{otherwise in $\elle{2}$},
\end{cases}
\end{equation}
where $I_{\Delta}$ is the indicator function of the convex set $\Delta$, defined in (\ref{indk}).

\begin{lemma}\label{l:basicp}
The functionals $F_p$, $p>3$, and $F$ are lower semicontinuous
in the strong $L^2$ topology.
Moreover, the functionals $F_p$ are equicoercive, i.e.\
there exists a constant $\alpha>0$ such that
\[
F_p(\vec{V})\geq \frac{\lambda}{2} \|\vec{V}\|_2^2-\alpha
\qquad
\forall \vec{V}\in \elle{2},\ \forall p>3\,.
\]
\end{lemma}

\begin{proof}
Let $3< p<\infty$, and let $(\vec{V}_k)\subset \elle{2}$
be a sequence converging to $\vec{V}$ in $L^2$. We have to prove that
\[
F_p(\vec{V})\leq \liminf_{k \to \infty}F_p(\vec{V}_k)\,.
\]
Without loss of generality we can assume that
$\vec{V}_k\in X_p$ for every $k$ and
\begin{equation}\label{f:limfin}
\lim_{k\to\infty} F_p(\vec{V}_k) = C < +\infty\,.
\end{equation}

We claim that $\vec{V}\in X_p$.
Namely, the divergence-free requirement is stable under
strong $L^2$ convergence.
Moreover, from (\ref{f:limfin}) and (\ref{robdd}) we have that
there exists a constant $C_1>0$ such that
$\|\curl\vec{V}_k\|_p \leq C_1$ for every $k\in\N$.
{}From Theorem~\ref{t:Xp}(i) we infer that the sequence
$(\vec{V}_k)$ is equicontinuous and equibounded in $\overline{A}$.
Hence we can pass to a subsequence, that we do not relabel,
converging uniformly to $\vec{V}$ in $\overline{A}$,
so that the boundary condition $\nt\times\vec{V} = 0$ on $\partial A$
holds.
Finally, from Remark~\ref{r:Xp},
$(\vec{V}_k)$ converges to $\vec{V}$ weakly in $\wp{p}$, and
$\|\curl\vec{V}\|_p \leq \liminf_k \|\curl\vec{V}_k\|_p$,
proving the claim.

By the convexity of $\gd$ and the convergence of
$\curl\vec{V}_k$ to $\curl\vec{V}$ in the weak $L^p$ topology,
we obtain that
\[
\int_A\gd(\curl\vec{V})^p \leq \liminf_{k\to \infty} \int_A\gd(\curl\vec{V}_k)^p,
\]
which implies  the semicontinuity inequality
$F_p(\vec{V})\leq C$.

The semicontinuity of $F$ can be proved in a similar way,
while the equicoerciveness of $(F_p)$ easily follows
choosing $\alpha = \frac{\lambda}{2}\|\vec{U}_0\|_2^2$.
\end{proof}

\begin{rk}\label{r:minFp}
For every $p>3$, the functional $F_p$ is lower semicontinuous in $\elle{2}$,
coercive,
and strictly convex in its effective domain $X_p$.
Hence it admits a unique minimizer $\vec{U}_p$,
which belongs to $X_p$.
Similarly, the functional $F$ admits a unique minimizer $\vec{U}\in X_2$
with $\gd(\curl\vec{U})\leq 1$.
\end{rk}

The mathematical justification of the Bad\'\i a and L\'opez model is
based on the fact that
the sequence $(\vec{U}_p)$ of
minimizers of $(F_p)$ converges to the minimizer $\vec{U}$ of $F$.
As is customary in material science,
we are going to prove this statement
using the technique of $\Gamma$-convergence (see \cite{Brai} and the
reference therein for examples of applications of the theory to
different models).

\begin{dhef}\label{d:gamma}
The sequence $(F_p)$ $\Gamma$-converges to the functional
$F$ (with respect to the strong $L^2$ topology) if the following
two conditions are satisfied.
\begin{itemize}
\item[(i)]
For every $\vec{V}\in \elle{2}$ there exists a sequence $(\vec{V}_p)$
(called a \textsl{recovering sequence}), converging to $\vec{V}$
in $\elle{2}$, and such that
\[
\limsup_{p\to\infty} F_p(\vec{V}_p) \leq F(\vec{V})\,.
\]
\item[(ii)]
For every sequence $(\vec{V}_p)\subset\elle{2}$ converging strongly to
$\vec{V}$, we have
\[
\liminf_{p\to\infty} F_p(\vec{V}_p) \geq F(\vec{V})\,.
\]
\end{itemize}
\end{dhef}

\begin{rk}\label{r:gc}
Since the functionals $F_p$ are equicoercive (see Lemma \ref{l:basicp}),
the convergence of the minimizers $(\vec{U}_p)$ to $\vec{U}$ will follow
if we prove that
$(F_p)$ $\Gamma$-converges to $F$ and
that $(\vec{U}_p)$ is a precompact sequence in $\elle{2}$
(see e.g.\ \cite{Brai}).
\end{rk}

The main tool needed in the proof of the $\Gamma$--convergence
of $(F_p)$ to $F$
is the following lemma.

\begin{lemma}\label{l:rholim}
Let $(\vec{V}_p)\in \elle{2}$ be a sequence converging to
$\vec{V}$ in $\elle{2}$. If
\[
\liminf_{p\to\infty} F_p(\vec{V}_p) < +\infty
\]
then
$\gd(\curl\vec{V})\leq 1$ a.e.\ in $A$.
\end{lemma}

\begin{proof}
Without loss of generality we can assume that
\[
\begin{split}
\liminf_{p\to\infty} F_p(\vec{V}_p) =
\lim_{p\to\infty} F_p(\vec{V}_p) = L < +\infty,\\
\curl\vec{V}_p\in \elle{p}\,,\quad F_p(\vec{V}_p)\leq L_1\qquad \forall p>3.
\end{split}
\]

We are going to prove that the sequence $(\vec{V}_p)$
is bounded in $\wp{q}$ for every $q>3$.
In view of Theorem~\ref{t:Xp}, it is enough show that
for every given $q>3$
there exists a constant $C>0$ such that
${\|\curl\vec{V}_p\|}_q \leq C$ for every $p>q$.

Given $3<q<p$, by (\ref{robdd}) and H\"older's inequality we have
\begin{equation}\label{f:gac1}
\begin{split}
\int_A |\curl\vec{V}_p|^q\, dx & \leq
\frac{1}{c_1^q}\int_{A} \gd(\curl\vec{V}_p)^q\, dx
\\ & \leq
\frac{1}{c_1^q}\left(\frac{p}{|A|}\right)^{q/p}|A|
\left(\int_A \frac{1}{p}\,\gd(\curl\vec{V}_p)^p\, dx\right)^{q/p}\,.
\end{split}
\end{equation}
On the other hand we get
\[
\int_A \frac{1}{p}\,\gd(\curl\vec{V}_p)^p\, dx \leq F_p(\vec{V}_p)
\leq L_1\,.
\]
Collecting all the previous estimates,
for every $p > q$ we obtain
\begin{equation}\label{f:gacn}
\int_A |\curl\vec{V}_p|^q\, dx \leq
\frac{ |A|}{c_1^q}\,
\left(\frac{p\, L_1}{|A|}\right)^{q/p}
\leq
\frac{|A|}{c_1^q}\,
\exp\left(\frac{q L_1}{e |A|}\right)\,.
\end{equation}
Then the sequence $(\vec{V}_p)$ is bounded in $\wp{q}$,
and hence $(\vec{V}_p)$ converges to $\vec{V}$ in the weak topology of
$\wp{q}$.

This fact and
the convexity of the function $\gd$ imply that
\begin{equation}\label{f:gacnpu}
\int_B \gd(\curl\vec{V})\, dx \leq \liminf_{p\to \infty}
\int_B \gd(\curl\vec{V}_p)\, dx
\end{equation}
for every open set $B\subseteq A$
(see e.g.\ \cite[Thm.~4.2.1]{Butt}).
On the other hand, by H\"older's inequality
\[
\begin{split}
\int_B \gd(\curl\vec{V}_p)\, dx & \leq
\left(\int_B \gd(\curl\vec{V}_p)^p\, dx\right)^{1/p}|B|^{(p-1)/p} \\
& \leq p^{1/p}\left(\frac{1}{p}\int_B \gd(\curl\vec{V}_p)^p dx\right)^{1/p}|B|^{(p-1)/p}\,.
\end{split}
\]
Then, using (\ref{f:gacn}) and (\ref{f:gacnpu}) we get
$\displaystyle\int_B \gd(\curl\vec{V})\, dx \leq |B|$.
Since $\gd(\curl\vec{V})\in L^1(A)$, we conclude that
\[
\gd(\curl\vec{V}(x))=\lim_{r\to 0^+} \frac{1}{|B_r(x)|}
\int_{B_r(x)} \gd(\curl\vec{V})\, dx \leq 1\,,
\]
for a.e.\ $x\in A$.
\end{proof}


\begin{theorem}\label{t:gammac}
The sequence $(F_p)$ $\Gamma$-converges to $F$.
\end{theorem}

\begin{proof}
(i)
Let $\vec{V}\in \elle{2}$.
It is not restrictive to assume that $F(\vec{V})<+\infty$,
so that $\vec{V}\in X_2$ and $\gd(\curl\vec{V})\leq 1$ a.e.
Then it is enough to choose $\vec{V}_p = \vec{V}$ as recovery sequence.

\noindent
(ii)
Let $\vec{V}_p \in \elle{2}$ be functions converging to
$\vec{V}$ in $\elle{2}$, such that
\[
\liminf\limits_{p\to\infty} F_p(\vec{V}_p)<+\infty.
\]
By Lemma~\ref{l:rholim} we have that $\gd(\curl\vec{V})\leq 1$
a.e.\ in $A$, and hence
\[
F(\vec{V}) = \int_A |\vec{V}-\vec{U}_0|^2\, dx
= \lim_{p\to\infty} \int_A |\vec{V}_p-\vec{U}_0|^2\, dx
\leq \liminf_{p\to\infty} F_p(\vec{V}_p)\,,
\]
concluding the proof.
\end{proof}
Notice that the $\Gamma$--convergence result is obtained under the sole
assumption that $\Uz\in\elle{2}$. Actually in our model
$\Uz=\vec{H}^n_i-\vec{H}_s(t^n_{i+1})$ so that $\Uz\in X_2+\vec{c}_0$,
$\cz\in\R^3$, and $\gd(\curl\Uz)\leq 1$. This additional regularity
is needed in order to recover the convergence of the minimizers.

\begin{theorem}
Let $\Uz\in \cz + X_2$, $\cz\in\R^3$, with $\gd(\curl\Uz)\leq 1$ a.e.\ in $A$.
Let $\vec{U}_p$ be the unique minimizer of $F_p$, $p>3$,
and let $\Usol$ be the unique minimizer of $F$.
Then $(\vec{U}_p)$ converges to $\vec{U}$ in the weak
topology of $\wp{q}$ for every $q>3$.
\end{theorem}

\begin{proof}
Reasoning as in (\ref{f:gac1}), for $3<q<p$ we have
\begin{equation}\label{f:gac2}
\int_A |\curl\vec{U}_p|^q\, dx \leq
\frac{1}{c_1^q}\left(\frac{p}{|A|}\right)^{q/p}|A|
\left(\int_A \frac{1}{p}\,\gd(\curl\vec{U}_p)^p\, dx\right)^{q/p}\,.
\end{equation}
On the other hand, by the minimality of $\vec{U}_p$ we get
\[
\int_A \frac{1}{p}\,\gd(\curl\vec{U}_p)^p\, dx \leq F_p(\vec{U}_p)
\leq F_p(\Uz-\cz)
\leq (1+\lambda\, |\cz|^2) |A|\,,
\]
so that
\[
\int_A |\curl\vec{U}_p|^q\, dx \leq
\frac{1}{c_1^q}\,
[p(1+\lambda\, |\cz|^2)]^{q/p} |A|
\leq
\frac{1}{c_1^q}\,|A|\,
\exp[q(1+\lambda\, |\cz|^2)/e]\,.
\]
Then the sequence $(\vec{U}_p)$ is bounded in
$\wp{q}$.
In particular,
it is precompact in $\elle{2}$.
Then, by Remark~\ref{r:gc} we obtain the
convergence of $(\vec{U}_p)$ to $\vec{U}$ in $\elle{2}$.
Finally, notice that the previous estimate implies that
$(\vec{U}_p)$ converges to $\vec{U}$ in
the weak topology of $W^{1,q}$ for every $q>3$.
\end{proof}

\section{Cylindrical superconductors}\label{s:bean2}

In what follows we shall consider a long cylindrical type--II superconductor
occupying the region $\Omega\times\R$, with a simply connected
cross section $\Omega\subset \R^2$, placed into a nonstationary uniform axial magnetic field
$\vec{H}_s(t)=(0,0,H_s(t))$.

Notice that, due to the parallel geometry, the problem admits a dimensional reduction.
Namely, we can assume that there exists a function
$u$ depending only on $x=(x_1,x_2)\in \Omega$ such that $\vec{H}=(0,0,H_s(t)+u(x,t))$, so that
$\curl \vec{H}=(u_{x_2}(x,t),-u_{x_1}(x,t), 0)$
and $\dive\vec{H}=0$.
In particular the current
density $\vec{J}$ is parallel to the cross section plane and its projection
on this plane is $(J_1,J_2)=(Du)^\bot$, where $(Du)^\bot$ is the
rotation of $-\frac{\pi}{2}$ of $Du$.
Hence, under the additional hypothesis that
the allowed region $\Delta$ is symmetric with respect to the $z=0$ plane, the constitutive law becomes
$Du\in K$, where $K$ is the rotation of the section $z=0$ of the set $\Delta$,
and $\vec{E}= (E_1,E_2, 0)$ with $(E_1,E_2)\in \partial I_K(Du)$.

In this case, if $c^n_i=H_s(t^n_i)$, and
$h^n_i(x)$ is the magnetic field intensity at time layer $i$, given by $H_s(t^n_i)+u(t^n_i,x)$,
the variational formulation is
\begin{equation}\label{varmod}
\min\left\{
\int_\Omega  (w-h^n_i)^2\, dx;\
w\in c^n_{i+1}+W^{1,1}_0(\Omega),\
Dw\in K\ \textrm{a.e.~in}\ \Omega
\right\}\,.
\end{equation}
In the following sections we shall give an explicit representation
of the unique minimizer $h^n_{i+1}$ and we shall show that the previous quasistatic evolution
converges to a function $h(x,t)$ (also explicit), which determines the evolution of the internal
magnetic field.

Further to the dimensional reduction and a rescaling, we obtain that
$h^n_{i+1}-c^n_{i+1}$ is the unique minimum point of the functional
\[
J(v)=\int_\Omega\left[
I_K(Dv)+
(v-\dato)^2\right]\, dx\,,\qquad  v\in W^{1,1}_0(\Omega),
\]
where $\dato=h^n_{i}-c^n_{i+1}$.

Although in this case we deal with an unbounded
superconductor $A=\Omega\times\R$,
the above reduction allows us to
prove, using the same arguments as in
Section~\ref{s:bean} (see also \cite{GNP}, Section 2),
a $\Gamma$--convergence result, and,
consequently, the convergence of the minimizers.
Even if the case $N=2$ is the
sole meaningful physical situation, the following results hold for
every dimension $N \geq 2$.

\begin{theorem}
Let $\Omega\subset\R^N$ be a bounded open set, let $K\subset\R^N$ be a nonempty, compact, convex set
with $0 \in \inte K$, and let $\gauge_K$ be the gauge function of $K$.
Let $G_p$, $p\geq 1$, and $G$ be the
functionals defined respectively by
\[
G_p(v)=
\begin{cases}
\displaystyle{\int_\Omega \frac{1}{p}\,\gauge_K(Dv)^p\, dx\,,} & v\in W^{1,p}_0(\Omega) \\
+\infty & \textrm{otherwise in }\ L^1(\Omega)
\end{cases}
\]
and
\[
G(v)=
\begin{cases}
0, & v\in W^{1,\infty}_0(\Omega),\ \gauge(Dv)\in K\\
+\infty & \textrm{otherwise in }\ L^1(\Omega)\,.
\end{cases}
\]
Then $(G_p)$ $\Gamma$--converges to $G$ in $L^r(\Omega)$ for every $r\geq 1$.
As a consequence, for a given
$\dato \in L^2(\Omega)$
the functionals
\[
{J}_p(v)=
\begin{cases}
\displaystyle{\int_\Omega \left[\frac{1}{p}\,\gauge_K(Dv)^p\, dx+
\left(v-\dato\right)^2\right]\, dx}\,, & v\in W^{1,p}_0(\Omega) \\
+\infty & \textrm{otherwise in }\ L^2(\Omega)
\end{cases}
\]
$\Gamma$--converge to the functional
\[
{J}(v)=
\begin{cases}
\displaystyle{\int_\Omega \left[
I_K(Dv)+
(v-\dato)^2\right]\, dx}\,, & v\in W^{1,1}_0(\Omega)\\
+\infty & \textrm{otherwise in }\ L^2(\Omega)\,,
\end{cases}
\]
in $L^2(\Omega)$.
If in addition  $ \dato \in c_0+W^{1,\infty}_0(\Omega)$, $c_0\in\R$, with $\gauge(D\dato)\leq 1$ a.e. in $\Omega$,
then the minimizers $u_p$ of $J_p$ converge to the unique minimizer $\uz$ of the functional $J$
in the weak topology of $W^{1,q}_0(\Omega)$ for every $q>1$.
\end{theorem}

\section{The Minkowski distance function}\label{s:distance}

Here and hereafter we shall deal with
\begin{equation}\label{f:Omega}
\Omega\subset\R^N\
\textrm{nonempty, bounded, open connected set
of class $C^2$},
\end{equation}
and
\begin{equation}\label{f:ipoK}
\begin{split}
& K\subset\R^N\,\ \textrm{nonempty, compact, convex set},\
0 \in \inte K, \\
&\partial K\ \textrm{of class}\ C^2\ \textrm{with strictly positive principal curvatures}\,.
\end{split}
\end{equation}

The polar set of $K$
is defined by
\[
K^0 = \{p\in\R^N;\ \pscal{p}{x}\leq 1\ \forall x\in K\}\,.
\]
We recall that, if $K$ satisfies (\ref{f:ipoK}), then
$K^0$ also satisfies (\ref{f:ipoK}), and $K^{00} = (K^0)^0 = K$
(see \cite[Thm.~1.6.1]{Sch}).

Since $K$ will be kept fixed, from now on we shall use for the gauge functions
the notation $\gauge = \gauge_K$ and $\pgauge=\gauge_{K^0}$

It can be readily seen that
the gauge function $\pgauge$ of the polar set $K^0$
coincides with the support function of the set $K$.
As a consequence, we have that
\begin{equation}\label{f:roroz}
\gauge(D\pgauge(\xi))=1\, \qquad
\forall\xi\in\R^N\setminus \{0\}\,.
\end{equation}


Let $\Lipr(\Omega)$ be the set of functions defined by
\[
\Lipr(\Omega) := \{u \in \Lip(\overline{\Omega}) \colon\ Du\in K
\ \textrm{a.e.\ in \ }\Omega\}.
\]
We recall that $u\in \Lipr(\Omega)$ if and only if
\begin{equation}\label{f:stimrho}
u(x)-u(y) \leq \pgauge(x-y)\,,\qquad
\end{equation}
for every $x$, $y\in\Omega$ joined by a segment contained in $\Omega$.

The main tools needed in the following sections are the
distances from the boundary of $\Omega$ associated to the
Minkowski structures induced by the gauge function of $K^0$
and $-K^0$ respectively.

\begin{dhef}\label{d:dist}
\ The Min\-kow\-ski distance from the
boundary of $\Omega$ is
\begin{equation}\label{f:d}
\dist(x) = 
\inf_{y\in\partial\Omega} \pgauge(x-y),
\qquad x\in\overline{\Omega}\,.
\end{equation}
Similarly, we define $\distm(x)$ as the distance from the boundary
induced by the gauge of $-K^0$.
\end{dhef}

Notice that the function $\distm$ coincides with $\dist$ only if
$K^0$ is symmetric with respect to the origin. In the remaining
part of this section we illustrate some features of the function
$\dist$. The analogous for $\distm$ can be obtained upon observing
that $-K^0$ is the polar set of $-K$.

Since $\partial\Omega$ is a compact subset of $\R^N$
and $\pgauge$ is a continuous function,
the infimum in the definition of $\dist$ is achieved.
We shall denote by $\proj(x)$ the set of projections of $x$
in $\partial\Omega$, that is
\begin{equation}\label{f:Pi}
\proj(x) = \{y\in\partial\Omega;\ \dist(x) = \pgauge(x-y)\},
\qquad x\in\overline{\Omega}.
\end{equation}

\begin{dhef}\label{d:sigma}
We say that $x\in\Omega$ is a
regular point of $\Omega$
if $\proj(x)$ is a singleton.
We say that $x\in\Omega$ is a
singular point of $\Omega$
if $x$ is not a regular point.
We denote by $\Sigma\subseteq\Omega$ the set
of all singular points of $\Omega$.
\end{dhef}

\begin{rk}\label{r:HJ}
It is well known that $\dist$ is a Lipschitz function
in $\Omega$, and that $\dist$ is differentiable
at $x\in\Omega$ if and only if $x$ is a regular point.
In addition, if $\proj(x) = \{y\}$ then
$D\dist(x)=D\pgauge(x-y)$, and hence, by (\ref{f:roroz}),
$D\dist\in K$ almost everywhere in $\Omega$
(see \cite{BaCD,CaSi,Li}; see also
Theorem~\ref{t:cm6}(i) below).
\end{rk}

\begin{rk}\label{r:maxm}
We recall that $\dist$ (resp. $\distm$) is the unique viscosity solution
of the Hamilton-Jacobi equation
$\gauge(Du) = 1$ (resp.\ $-\gauge(Du) =-1$) in $\Omega$,
with boundary condition
$u=0$ on $\partial\Omega$.
As a consequence of the maximality property of the viscosity
solutions, we have that
\[
-\distm\leq u \leq \dist\,, \qquad \forall u\in\Lipr(\Omega),\ u=0\ \textrm{on}\ \partial\Omega
\]
(see \cite{Li}). Moreover it can be easily checked that
\[
\begin{aligned}
u  &\geq -\distm\,, &\quad \forall u\in\Lipr(\Omega),\ u\geq 0\ \textrm{on}\ \partial \Omega\,, \\
u  &\leq \dist\,,  & \quad \forall u\in\Lipr(\Omega),\ u\leq 0\ \textrm{on}\ \partial \Omega\,,
\end{aligned}
\]
(see e.g.\ \cite{BaCD}, Theorem 5.9).
\end{rk}

In the following theorem
we collect
all the results proved in \cite{CMf} that are
relevant for the subsequent analysis.

\begin{theorem}\label{t:cm6}
Assume that $\Omega$ and $K$ satisfy respectively $(\ref{f:Omega})$
and $(\ref{f:ipoK})$. Then the following hold.
\begin{itemize}
\item[(i)] $\overline{\Sigma}\subset\Omega$,
and the Lebesgue measure of
$\overline{\Sigma}$ is zero.

\item[(ii)] 
Let $x\in\Omega$ and $\macroy \in\proj(x)$.
Then $\dist$ is differentiable at every $z$
along the segment jointing $y$ to $x$ (without endpoints),
and
$
D\dist(z) = \frac{\nor(\macroy )}{\gauge(\nor(\macroy ))}
$,
where $\nor(\macroy )$ is the (Euclidean) inward normal unit vector
to $\partial\Omega$ at $\macroy $.

\item[(iii)] 
The function $\dist$ is of class $C^2$ in
$\overline{\Omega}\setminus\overline{\Sigma}$.
\end{itemize}
\end{theorem}

\begin{proof}
See  Remark~4.16, Corollary~6.9,
Lemma~4.3
and Theorem~6.10 in \cite{CMf}.
\end{proof}

{}At any point $\macroy \in\partial\Omega$
there is a unique inward ``normal'' direction $p(\macroy )$
with the properties
$\proj(\macroy +t p(\macroy )) = \{\macroy \}$
and $\dist(\macroy +t p(\macroy )) = t$
for $t\geq 0$ small enough
(see \cite[Remark~4.5]{CMf}).
More precisely, these properties
hold true for $p(\macroy ) = D\gauge(\nor(\macroy ))$
and for every $t\in [0, \len(\macroy ))$,
where $\len(\macroy )$ is defined by
\begin{equation}\label{f:len}
\len(\macroy ) =
\min\{t\geq 0;\
\macroy  + t D\gauge(D\dist(x))\in\overline{\Sigma}\}
\end{equation}
(see
\cite[Propositions~4.4 and~4.8]{CMf} and \cite[Lemma~2.2]{LN}).
{}From Theorem~\ref{t:cm6}(ii) and the positive $0$-homogeneity
of $D\gauge$ it is plain that
$D\gauge(\nor(\macroy )) = D\gauge(D\dist(\macroy ))$.

{}From Theorem~\ref{t:cm6}(iii), the function $\dist$
is of class $C^2$ on
$\partial\Omega$.
Then we can define the function
\begin{equation}\label{f:W}
W(\macroy) = -D^2\gauge(D\dist(\macroy))\, D^2\dist(\macroy)\,,\qquad
\macroy\in
\partial\Omega\,.
\end{equation}
For any $\macroy \in\partial\Omega$ let $T_{\macroy }$ denote
the tangent space to $\partial\Omega$ at $\macroy $.
It can be proved that for every $v\in T_y$,
one has $W(y)\, v\in T_y$.
Hence, we can define the map
\begin{equation}\label{f:bw}
\bw(y)\colon T_y\to T_y,\quad
\bw(y)\, w = W(y)\, w,
\end{equation}
that can be identified with a linear application
from $\R^{N-1}$ to $\R^{N-1}$.

Although
the matrix $\bw(\macroy)$ is not in general symmetric,
its eigenvalues
are real numbers,
and so its eigenvectors are real
(see \cite{CMf}, Remark~5.3).
The eigenvalues of $\bw(\macroy)$
have an important geometric interpretation.

\begin{dhef}\label{d:curv}
Let $\macroy \in\partial\Omega$.
The anisotropic curvatures of $\partial\Omega$ at $\macroy $,
with respect to the Minkowski norm $\dist$,
are the eigenvalues
$\curvg_1(\macroy )\leq\cdots\leq\curvg_{n-1}(\macroy )$
of $\bw(\macroy )$.
\end{dhef}

In what follows we shall denote by $\Pi^-$, $\curvg^-_j$, $j=1,\ldots,n-1$, $\Sigma^-$ and
$\len^-$,
respectively
the projection, the anisotropic curvatures and the singular set
associated to $\distm$, and the normal distance to the cut locus $\Sigma^-$.

\section{A minimum problem}\label{s:minp}


Fixed $\dato\in \Lipr(\Omega)$, we want to minimize the functional
\begin{equation}\label{jfunct}
J(v)=\int_\Omega [I_K(Dv)+(v-\dato)^2]\, dx\,, \quad v\in W^{1,1}_0(\Omega)\,.
\end{equation}
Let us consider the partition
$\Omega = \Omega^+\cup\Omega^-\cup \Omega^0$, with
\[
\begin{split}
\Omega^0  =\{-\distm \leq \dato \leq \dist\}\,,\
\Omega^+ =\{\dato > \dist\}\,, \
\Omega^-=\{ \dato <- \distm\}\,,
\end{split}
\]
and define the function
\begin{equation}\label{umin}
u(x)=
\begin{cases}
\dist(x) & x\in \Omega^+ \\
-\distm(x) & x\in \Omega^- \\
\dato(x) & x \in \Omega^0\,.
\end{cases}
\end{equation}
The aim of this section is to show that $u$ is the unique minimizer
of $J$ in $W^{1,1}_0(\Omega)$.

\begin{rk}\label{r:noint}
If $x \in \overline{\Omega}^+ \cap \overline{\Omega}^-$, then
$\dato(x)=\dist(x)=-\distm(x)$ and hence $x\in \partial \Omega$
and $\dato(x)=0$. As a consequence, we have
$u(x)=\min((\max(\dato(x), -\distm(x)), \dist(x)))$.
This implies that
$u\in \Lipr(\Omega)$, $u=0$ on $\partial\Omega$,
$\gauge(Du)=1$ a.e.\ in $\Omega^+\cup\Omega^-$, and $u$ is a viscosity
solution of $\gauge(Du)=1$ in $\Omega^+$ and of $-\gauge(Du)=-1$ in $\Omega^-$.
\end{rk}

Let us define the following two subsets of $\partial\Omega$:
\begin{equation}\label{f:gpm}
\Gamma^+=\{y\in\partial\Omega \colon \dato(y)>0\},\qquad
\Gamma^-=\{y\in\partial\Omega \colon \dato(y)<0\}\,,
\end{equation}
and the functions
\begin{equation}\label{f:lpm}
\begin{split}
\lambda(y) & = \sup\{s\in [0,\len(y))\colon \dato(y+tD\gauge(\nor(y)))>t\
\forall\ t\in [0,s)\}\,,\\
\lambda^-(y) & = \sup\{s\in [0,\len^-(y))\colon \dato(y+tD\gaugem(\nor(y)))<-t\
\forall\ t\in [0,s)\}\,.
\end{split}
\end{equation}
Recall that $\len(y)$ and $\len^-(y)$
are the normal distance to the cut locus,
corresponding to $\dist$ and $\distm$ respectively, that is the length of the ``normal'' ray
starting from $y$ and ending on $\overline{\Sigma}$ (see (\ref{f:len})).

\begin{rk}\label{r:lambda}
It can be easily checked that $\lambda$ and $\lambda^-$ are lower semicontinuous function
in $\partial \Omega$. In general we cannot expect the continuity of these functions, as it
is shown by the following example.
Let $\Omega=B_1(0)\subseteq \R^2$, $K=\overline{B}_1(0)\subseteq \R^2$, and
define
$\dato(r \cos \theta,r \sin \theta)=1-r \cos \theta$
for $0\leq r<1$, $\theta\in [0,2\pi)$.
We have
$|D\dato(r \cos\theta, r\sin\theta)|^2=\cos^2 \theta + \sin^2 \theta=1$
($0<r<1$), so that $\dato\in \Lipr(\Omega)$.
Moreover, if $y=(\cos\theta, \sin\theta)\in \partial \Omega$, $\theta\in [0,2\pi)$,
$\lambda$ is
the lower semicontinuous function given by
\[
\lambda(y)=
\begin{cases}
0\,, & \theta=0\,, \\
1\,, & \theta\in (0,2\pi)\,.
\end{cases}
\]
\end{rk}

\smallskip

The following simple lemma provides a
characterization of $\Omega^+$ and $\Omega^-$ in terms of the rays starting from points
of $\Gamma^+$ and $\Gamma^-$ respectively.

\begin{lemma}\label{l:opm}
Let $y\in\Gamma^+$ be such that
$\lambda(y)<\len(y)$.
Then $y+t D\gauge(\nor(y))\in\Omega^0$
for every $t\in [\lambda(y),\len(y)]$.
Analogously,
let $y\in\Gamma^-$ be such that
$\lambda^-(y)<\len^-(y)$.
Then
$y+t D\gaugem(\nor(y))\in\Omega^0$
for every $t\in [\lambda^-(y),\len^-(y)]$. As a consequence, we have
\begin{equation}\label{f:splom}
\begin{split}
\Omega^+\setminus \overline{\Sigma}^+ & =\{y+t D\gauge(\nor(y))\colon y\in\Gamma^+,\
t\in (0, \lambda(y))\},\\
\Omega^-\setminus \overline{\Sigma}^- & =\{y+t D\gaugem(\nor(y))\colon y\in\Gamma^-,\
t\in (0, \lambda^-(y))\}.
\end{split}
\end{equation}
\end{lemma}

\begin{proof}
Let $y\in\Gamma^+$ and assume that $\lambda(y)<\len(y)$.
{}From the very definition of $\lambda(y)$ and the continuity of $\dato$, we have that
$\dato(y+\lambda(y)\, D\gauge(\nor(y))) = \lambda(y)$.
For every $t\in [\lambda(y), \len(y)]$ we have that
\[
\dist(y+tD\gauge(\nor(y)))= t = \dato(y+\lambda(y)\, D\gauge(\nor(y))) + t-\lambda(y) \geq
\dato(y+t\, D\gauge(\nor(y)))\,,
\]
where the last inequality follows from the assumption $\dato\in\Lip_\rho(\Omega)$
(see (\ref{f:stimrho})), and from the fact that $\pgauge((t-\lambda)D\gauge(\nor(y)))= t-\lambda$
(see (\ref{f:roroz})).
The inequality
$\dato(y+t\, D\gauge(\nor(y))) \geq -\distm(y+t\, D\gauge(\nor(y)))$
follows easily from the fact that
$\dato(y) > 0$ and $\dato\in\Lipr(\Omega)$ (see Remark \ref{r:maxm}).
The case $y\in\Gamma^-$ can be handled in a similar way.
\end{proof}


The following change of variable formula can be proved as in
\cite[Theorem~7.1]{CMf}.

\begin{theorem}\label{t:chvar}
Let $\Phi, \Phi^-\colon\partial\Omega\times\R\to\R^N$
be the maps defined respectively by
\[
\Phi(y,t) = y + t\, D\gauge(\nor(y)), \quad
\Phi^-(y,t) = y + t\, D\gaugem(\nor(y)),
\ \quad (y,t)\in\partial\Omega\times\R.
\]
Then
\[
\begin{split}
\int_{\Omega^+} h(x)\, dx
&=
\int_{\Gamma^+} \gauge(\nor(y))\, \left[
\int_{0}^{\lambda(y)}
h(\Phi(y,t))\,
\prod_{i=1}^{n-1} (1-t\, \curvg_i(y))
\, dt\right]\,d\haus(y)
\\
\int_{\Omega^-} h(x)\, dx
&=
\int_{\Gamma^-} \gaugem(\nor(y))\, \left[
\int_{0}^{\lambda^-(y)}
h(\Phi^-(y,t))\,
\prod_{i=1}^{n-1} (1-t\, \curvg^-_i(y))
\, dt\right]\,d\haus(y)
\end{split}
\]
for every $h\in L^1(\Omega)$.
\end{theorem}

In what follows we will set
\[
M_{x}(t) = \prod_{i=1}^{n-1}\frac{1-t\curvg_i(x)}{1-\dist(x)\curvg_i(x)}\,, \quad
M_{x}^-(t) = \prod_{i=1}^{n-1}\frac{1-t\curvgm_i(x)}{1-\distm(x)\curvgm_i(x)}\,.
\]
{}From Lemma~7.3 in \cite{CMf} we deduce that there exists a positive constant $M_0$
such that
\begin{equation}\label{f:estiM}
0\leq M_x(t)\leq M_0\quad
\forall x\in\Omega\setminus\overline{\Sigma},\ \dist(x)\leq t < \len(y),
\end{equation}
where $\proj(x) = \{y\}$.
An analogous estimate holds for $M_x^-(t)$.

Using Theorem \ref{t:chvar} and (\ref{f:splom}), we are able to
show that the function $u$ defined in (\ref{umin}), coupled with
an explicit function $\vf$, solves a system of PDEs of Monge--Kantorowich
type. \textit{A posteriori} this system will be understood as a
first order necessary condition fulfilled by the minimizer of $J$.

\begin{lemma}\label{l:eul}
Let $u\in \Lipr(\Omega)$ be the function defined in (\ref{umin}),
and let $\vf\colon\Omega\to\R$ be the function defined by
\begin{equation}\label{f:defv}
\vf(x)=
\begin{cases}
\displaystyle
\int_{\dist(x)}^{\lambda(\macroy )} [\dato(\macroy +t\,D\gauge(D\dist(x)))-t]\,
M_{x}(t) dt\,, & x\in\Omega^+\setminus \overline{\Sigma}^+,\ \Pi(x)=\{y\}\\ \\
\displaystyle
-\int_{\distm(x)}^{\lambda^-(\macroy )} [\dato(\macroy +t\,D\gaugem(D\dist(x)))+t]\,
M_{x}^-(t)
\, dt\,, & x\in\Omega^-\setminus \overline{\Sigma}^-,\ \Pi^-(x)=\{y\}\\
0 & \textrm{otherwise}\,.
\end{cases}
\end{equation}
Then $\vf\in C_b(\Omega)$, $\vf\geq 0$, and the pair $(u,\vf)$ solves the system of PDEs
\begin{equation}\label{f:MK}
\begin{cases}
-\dive(\vf\, D\gauge(Du)) = \dato-u
&\textrm{in $\Omega$},\quad \textrm{(distributional)}\\
\gauge(Du)\leq 1
&\textrm{a.e.\ in $\Omega$},\\
\gauge(Du) = 1
&\textrm{a.e.\ in $\{\vf>0\}$},\\
u=0 & \textrm{in $\partial\Omega$}\,.
\end{cases}
\end{equation}
\end{lemma}

\begin{proof}
Let us extend the functions $\lambda$ and $\lambda^-$
to $\Omega^+\setminus\Sigma^+$ and $\Omega^-\setminus\Sigma^-$
respectively by setting
$\lambda(x):= \lambda(y)$ when $\proj(x) = \{y\}$ and
$\lambda^-(x):= \lambda^-(y)$ when $\proj^-(x) = \{y\}$.
Since $\lambda(x)\geq\dist(x)$ for every $x\in\Omega^+\setminus\Sigma^+$
and $\lambda^-(x)\geq\distm(x)$ for every $x\in\Omega^-\setminus\Sigma^-$,
it is plain that $\vf\geq 0$.

Let us prove that $\vf$ is a continuous function.
It is convenient to rewrite $\vf$ in $\Omega^+\setminus \overline{\Sigma}^+$ in the following way:
\begin{equation}\label{f:defvp}
\vf(x)=
\int_{\dist(x)}^{\len(\macroy )}
[t\vee\dato(\macroy +t\,D\gauge(D\dist(x)))-t]\,
M_{x}(t) dt\,, \quad x\in\Omega^+\setminus \overline{\Sigma}^+,\ \Pi(x)=\{y\}.
\end{equation}
Since the maps $\dist$, $D\dist$, $\dato$, $\len$, and $x\mapsto\proj(x)$ are continuous
in $\Omega^+\setminus \overline{\Sigma}^+$,
it follows that also $\vf$ is continuous in $\Omega^+\setminus \overline{\Sigma}^+$.

If $x_0\in\Omega^+\cap\overline{\Sigma}^+$,
then $\lambda(x_0) = \len(x_0) = \dist(x_0)$,
and the map $\lambda$ is continuous at $x_0$.
Namely, since $x_0\in\overline{\Sigma}^+$, we have that
$\len(x_0) = \dist(x_0)$.
Moreover, since $x_0\in\Omega^+$, we have that
$\dato(x_0) > \dist(x_0)$, so that $\dato > \dist$ on every segment
$[y,x_0]\subset\overline{\Omega}$ with $y\in\partial\Omega$.
Thus $\lambda(x_0)\geq\dist(x_0)$.
Finally, from the lower semicontinuity of $\lambda$
and the fact that $\lambda(x)\leq\len(x)$ for every $x\in\Omega$,
we conclude that
\[
\dist(x_0)\leq\lambda(x_0)\leq\liminf_{x\to x_0}\lambda(x)\leq
\limsup_{x\to x_0}\lambda(x)\leq
\lim_{x\to x_0}\len(x) = \len(x_0) = \dist(x_0),
\]
so that
$\lambda(x_0) = \dist(x_0) = \lim_{x\to x_0}\lambda(x)$.

{}From the very definition of $\vf$ and
(\ref{f:estiM})
we have that
\[
0\leq \vf(x)\leq C[\lambda(x)-\dist(x)]\qquad
\forall x\in \Omega^+\setminus\overline{\Sigma}^+.
\]
We remark that the same estimate trivially holds also on
$\Omega^+\cap\overline{\Sigma}^+$ since $\vf$ vanishes.
If $x_0\in\Omega^+\cap\overline{\Sigma}^+$,
we have that
\[
\limsup_{x\to x_0}\vf(x)\leq
\limsup_{x\to x_0} C[\lambda(x)-\dist(x)] = 0,
\]
hence $\lim_{x\to x_0}\vf(x) = 0 = \vf(x_0)$.

Up to now we have proved that $\vf$ is continuous in $\Omega^+$.
A similar argument shows that $\vf$ is continuous in $\Omega^-$.
Since $\vf$ vanishes on $\Omega^0$ and
$\partial\Omega^+\cap\partial\Omega^-\cap\Omega=\emptyset$
(see Remark~\ref{r:noint}),
it remains to prove that $\vf$ is continuous on
$\Omega\cap\partial\Omega^+$ and $\Omega\cap\partial\Omega^-$.

Let $x_0\in\Omega\cap\partial\Omega^+$.
Since $x_0\in \Omega^0$, we have that $v(x_0)=0$.
Moreover, by definition of $\Omega^+$,
$\dato(x_0) = \dist(x_0)$.
Since $\dato\in\Lipr(\Omega)$,
by (\ref{f:stimrho}) and (\ref{robdd}),
given $\epsilon > 0$ we have that
\[
\dato(x) \leq \dist(x) + 2 c_2 \epsilon
\qquad\forall x\in\Omega\cap B_{\epsilon}(x_0).
\]
Moreover, if $x\in(\Omega^+\cap B_{\epsilon}(x_0))\setminus\overline{\Sigma}^+$,
we have that $\dato\leq\dist + 2c_2\epsilon$
along the segment joining $x$ to its cut point.
As a consequence,
the integrand in (\ref{f:defvp}) is bounded from above by
$2 c_2\epsilon\, M_0$, so that
\[
0\leq\vf(x)\leq C \epsilon
\qquad\forall x\in(\Omega^+\cap B_{\epsilon}(x_0))\setminus\overline{\Sigma}^+.
\]
Since $\vf = 0$ in $\Omega^0\cup\overline{\Sigma}^+$,
this inequality clearly implies that
\[
\lim_{x\to x_0}\vf(x) = 0 = \vf(x_0),
\]
i.e.\ $\vf$ is continuous at $x_0$.
The continuity in $\Omega\cap\partial\Omega^-$
can be proved in a similar way.

\smallskip
Let us prove that the pair $(u,\vf)$ is a solution to (\ref{f:MK}).
Upon observing that $\{\vf>0\}=\Omega^+ \cup \Omega^-$,
by Remark~\ref{r:noint}
we have only to
prove that $(u,\vf)$ solves
\[
\int_\Omega \vf\pscal{D\gauge(Du)}{D\varphi}\, dx = \int_\Omega (\dato-u)\varphi\, dx\,,
\]
for all $\varphi\in C^{\infty}_0(\Omega)$.
We have that
\[
\int_\Omega (\dato-u)\varphi\, dx=
\int_{\Omega^+} (\dato-u)\varphi\, dx+\int_{\Omega^-} (\dato-u)\varphi\, dx\,.
\]
By using the change of variables stated in Theorem~\ref{t:chvar}, we get
\[
\begin{split}
\int_{\Omega^-} (\dato-u)\varphi =
\int_{\Gamma^-}\gaugem(\nor)
\left[\int_0^{\lambda^-}(\dato(\Phi^-)-u(\Phi^-))\varphi(\Phi^-)
\prod_{i=1}^{n-1} (1-t\, \curvg_i^-)
\, dt\right]\,d\haus\,.
\end{split}
\]
An integration by parts leads to
\[
\begin{split}
&\int_{0}^{\lambda^-}
(\dato(\Phi^-)-u(\Phi^-))\,\varphi(\Phi^-)\,
\prod_{i=1}^{n-1} (1-t\, \curvg_i^-)
\, dt  \\
& =
\int_{0}^{\lambda^-}
\pscal{D\varphi(\Phi^-)}{D\gaugem(\nor)}
\left[\int_{t}^{\lambda^-}
(\dato(\Phi^-)-u(\Phi^-))\,
\prod_{i=1}^{n-1} (1-s\, \curvg_i^-)\, ds\right]
\, dt\,.
\end{split}
\]
Since $u(\Phi^-(y,s)) = -s$ for every $s\in [0, \lambda^-(y)]$,
we have that
\[
-\vf(\Phi^-(y,t)) =
\int_{t}^{\lambda^-(y)}
(\dato(\Phi^-(y,s))-u(\Phi^-(y,s)))\,
\prod_{i=1}^{n-1}\frac{1-s\, \curvg_i^-(y)}%
{1-t\, \curvg_i^-(y)}\, ds
\]
for every $y\in\Gamma^-$ and $t\in [0,\lambda^-(y))$.
Then we obtain
\[
\begin{split}
\int_{\Omega^-} & (\dato-u)\varphi\, dx
=
-\int_{\Gamma^-}\gaugem(\nor)\left[
\int_{0}^{\lambda^-}
\pscal{D\varphi(\Phi^-)}{D\gaugem(\nor)}
\vf(\Phi)
\prod_{i=1}^{n-1} (1-t\, \curvg_i^-)\, dt\right]\, d\haus\,.
\end{split}
\]
Finally, again from Theorem~\ref{t:chvar} and the fact that
\[
D\gaugem(\nor(y))=
- D \gauge (-\nor (y))=-D \gauge(-D\distm(y))=-D\gauge(Du(y))\,,
\]
we get
\[
\int_{\Omega^-} (\dato-u)\varphi\, dx =
\int_{\Omega^-}
\vf\, \pscal{D\gauge(Du)}{D\varphi}\, dx\,,
\]
for every $\varphi\in C^{\infty}_c(\Omega)$.
In the same way,
using the change of variables in $\Omega^+$,
we get
\[
\int_{\Omega^+} (\dato-u)\varphi =
\int_{\Omega^+}
\vf\, \pscal{D\gauge(Du)}{D\varphi}\,,
\]
for every $\varphi\in C^{\infty}_c(\Omega)$,
concluding the proof.
\end{proof}

Now the absolute minimality of $u$ follows easily.
\begin{theorem}\label{t:exmin}
The function $u$ defined in (\ref{umin}) is the unique
minimizer of $J$ in $W^{1,1}_0(\Omega)$.
\end{theorem}

\begin{proof}
Let $\vf$ be the function defined in (\ref{f:defv}). Since $\vf\geq 0$,
for almost every $x\in \Omega^+\cup\Omega^-$
the vector $2 \vf(x) D\gauge(Du(x))$ belongs to the subdifferential
of $I_K$ at $Du(x)$ (see (\ref{subindk})).
On the other hand, $\vf = 0$ in $\Omega^0$.
Then we have
\[
I_K(Dw(x))-I_K(Du(x)) \geq 2 \vf(x) \pscal{D\gauge(Du(x))}{Dw(x)-Du(x)}
\ \textrm{a.e.\ in\ }\Omega\,,
\]
for every $w\in W^{1,1}_0(\Omega)$.
Hence, using Lemma \ref{l:eul} we obtain that for every
$w\in W^{1,1}_0(\Omega)$
\begin{equation}\label{f:minprf}
\begin{split}
J(w) - J(u)
& \geq {}
\int_\Omega 2 \vf \pscal{D\gauge(Du)}{Dw-Du} +
\int_\Omega (w-\dato)^2
+ \int_\Omega (u-\dato)^2
\\ = {} &
\int_\Omega [2(w-u)(\dato-u)+(w-\dato)^2-(u-\dato)^2]
\\ = {} &
\int_\Omega(u-w)^2\,,
\end{split}
\end{equation}
which implies that $u$ is the unique minimizer of $J$ in $W^{1,1}_0(\Omega)$.
\end{proof}

\section{The quasistatic evolution}\label{s:evol}
Let $T>0$ be fixed. For every $n\in\N^+$, we take a partition
$P^n=\left\{t^n_i\right\}_{i=0}^{k(n)}$
of the interval $[0,T]$, and we set $\delta t^n_i=t^n_i-t^n_{i-1}$.

The starting configuration of the superconductor at time $t=0$ is given by
$h(x,0)=h_0(x)\in \Lipr(\Omega)$, with $h_0=H_s(0)$ on $\partial \Omega$. 
Using the variational formulation (\ref{varmod}) and Theorem \ref{t:exmin}, we have that
the internal magnetic field $h^n_{i+1}(x)=h(x, t^n_{i+1})$ is given by
\begin{equation}\label{f:akkani}
h^n_{i+1}(x)=\min[\max[h^n_i(x), c^n_{i+1}-\distm(x)], c^n_{i+1}+\dist(x)]\,.
\end{equation}
Since the model is valid only when the external field varies slowly in time,
it is not restrictive to assume that $H_s$ is a continuous function obtained
by the concatenation of a finite number of monotone functions of class $C^1$.
Hence from now on we shall consider only the case
\begin{equation}\label{f:HS}
\textrm{$H_s$ is a monotone
function of class $C^1([0,T])$.}
\end{equation}
In this case formula (\ref{f:akkani}) simplifies.
Namely, if we assume that $H_s$ is a nondecreasing
function we have that $c^n_i \leq c^n_{i+1}$.
In particular, at the first
time step we have
\[
h^n_1(x)=\min[\max[h_0(x), c^n_{1}-\distm(x)], c^n_{1}+\dist(x)]\,,
\]
and $h_0(x)-c^n_{1} \in \Lipr(\Omega)$ with
$h_0(x)-c^n_{1}=c^n_0-c^n_{1}\leq 0$ on $\partial \Omega$.
Hence $h_0(x)-c^n_{1} \leq \dist (x)$ in $\Omega$ (see Remark \ref{r:maxm}), which
implies
\[
h^n_1(x)=\max[h_0(x), c^n_{1}-\distm(x)]\,.
\]
At the second time step we get
\[
h^n_2(x)=\min[\max[h^n_1(x), c^n_{2}-\distm(x)], c^n_{2}+\dist(x)]\,.
\]
Since
$h^n_1(x)-c^n_{2}=\max[h_0(x)-c^n_{2}, -\distm(x)+c^n_1-c^n_2]$, with
$c^n_1-c^n_2\leq 0$, we have that
\[
h^n_2(x)=\max[h_0(x), c^n_{2}-\distm(x)]\,.
\]
In general we get
\begin{equation}\label{f:akkanicres}
h^n_{i}(x)=\max[h_0(x), c^n_{i}-\distm(x)]\,.\quad \textrm{(nondecreasing external field)}
\end{equation}
With a similar argument, we obtain that if $H_s$ is nonincreasing in time, the
$i$--th time step leads to
\begin{equation}\label{f:akkanidecre}
h^n_{i}(x)=\min[h_0(x), c^n_{i}+\dist(x)]\,.\quad \textrm{(nonincreasing external field)}
\end{equation}
Now we easily obtain the evolution of the internal magnetic field, taking the
limit as $n$ goes to $+\infty$ of the quasistatic approximation
\[
h^n(x,t)=h^n_i(x)\,,\ t\in [t^n_i, t^n_{i+1})\,.
\]

\begin{theorem}\label{t:limitsol}
If $H_s$ is a nondecreasing function in $[0,T]$, then
the sequence
$(h^n(x,t))$ converges uniformly in $\overline{\Omega} \times [0,T]$ to
\[
h(x,t)=\max[h_0(x), -\distm(x)+H_s(t)]\,.
\]
If $H_s$ is a nonincreasing function in $[0,T]$, then
the sequence
$(h^n(x,t))$ converges uniformly in $\overline{\Omega}\times [0,T]$ to
\[
h(x,t)=\min[h_0(x), \dist(x)+H_s(t)]\,.
\]
\end{theorem}

\begin{proof}
Assume that $H_s$ is nondecreasing (the other case being similar). Given $t\in [0,T]$, and $n\in\N$,
we have that $t\in [t^n_i, t^n_{i+1})$, and
\[
0\leq h(x,t)-h^n(x,t)\leq H_s(t)-H_s(t^n_i) \leq H_s(t^n_{i+1})-H_s(t^n_i)
\leq \delta t^n \max_{t\in [0,T]}H_s'(t)\,,
\]
where $\delta t^n$ is the size of the partition $P^n$. Hence the uniform
convergence in $\overline{\Omega}\times [0,T]$ follows.
\end{proof}


\begin{rk}
Theorem \ref{t:limitsol} shows that the evolution obtained
by the direct optimization method at discrete times proposed by
Bad\'\i a and L\'opez in \cite{BLd} coincides with the one obtained in
\cite{BaPrb,BaPr} by
Barrett and Prigozhin as a solution of an evolutionary
variational inequality.
\end{rk}

For what concerns the electric field induced in the superconductor, we already
know that $\vec{E}=(E_1,E_2,0)$ with $(E_1,E_2)\in \partial I_K(Dh)$. This means that
$\vec{E}=0$ whenever $Dh \in \textrm{int}K$, and there exists
$w(x,t)\geq 0$ such that $(E_1,E_2)=w(x,t) D\gauge(Dh(x,t))$
for almost every $x\in \Omega$ and for every $t\in [0,T]$.
Our approach allows us to compute explicitly the quasistatic evolution
of the dissipated power $w(x,t)$ as the dual function appearing in the necessary
conditions stated in Lemma~\ref{l:eul}.

Namely, from Lemma~\ref{l:eul}, for every $n\in\N$
and for every $i=1, \ldots, k(n)$, there exists $v^n_{i}\in C_b(\Omega)$
such that
\begin{equation}\label{f:faradisc}
-\dive\left(\frac{v^n_{i}}{\dt}\,D\gauge(Dh^n_{i})\right)=
\frac{h^n_{i-1}-h^n_{i}}{\dt}\,,
\end{equation}
in the sense of distributions in $\Omega$.
For our convenience, we also set $v^n_0 = 0$.
The function $v^n_i$ is unique, as stated in the following result.

\begin{prop}\label{p:vni}
Let $v \in C_b(\Omega)$ be a function
such that
\begin{equation}\label{f:vvv}
\begin{cases}
-\dive\left(v\,D\gauge(Dh^n_{i})\right)=
{h^n_{i-1}-h^n_{i}}\,,
&\textrm{in $\Omega$ (distributional)} ,\\
\gauge(D h^n_i) = 1
&\textrm{a.e.\ in $\{v > 0\}$}.
\end{cases}
\end{equation}
If $c^n_i\geq c^n_{i-1}$, then $v=v^n_i$, where
\[
v^n_i(x)=
\begin{cases}
\displaystyle\int_{\distm(x)}^{\len^-(y)}{(h^n_i-
h^n_{i-1})(y+s D\gaugem(\nor(y)))}
M_{x}^-(s)\, ds\,, & x\in \Omega \setminus \overline{\Sigma}^-\,,\ \proj^-(x)=\{y\},\\
0 & x\in \overline{\Sigma}^-
\end{cases}
\]
whereas, if $c^n_i\leq c^n_{i-1}$, then $v=v^n_i$, where
\[
v^n_i(x)=
\begin{cases}
\displaystyle\int_{\dist(x)}^{\len(y)}{(h^n_{i-1}-
h^n_{i})(y+s D\gauge(\nor(y)))}
M_{x}(s)\, ds\,, & x\in \Omega \setminus \overline{\Sigma}\,,\ \proj(x)=\{y\},\\
0 & x\in \overline{\Sigma}\,.
\end{cases}
\]
\end{prop}

\begin{proof}
Assume that $c^n_{i-1}-c^n_i\geq 0$ (the other case being similar),
and let $v\in C_b(\Omega)$ be a solution of (\ref{f:vvv}).
By construction we have that the function
$f := h^n_{i-1}-h^n_{i}$ is bounded, continuous and nonnegative in $\Omega$,
and $h^n_i = \dist$ in $\spt(f)$.
Using the same arguments of Lemma~6.3 in \cite{CMg}, we can prove that
the pair $(\dist, v)$ satisfies
\[
-\dive\left(v\,D\gauge(D\dist)\right)=f,
\]
in the sense of distributions in $\Omega$.
Now we can apply Propositions~6.5 and 6.7 in \cite{CMg} in order to
conclude that $v = v^n_i$.
\end{proof}

Formula (\ref{f:faradisc})
can be understood as a discretized version of Faraday's law, so that
the functions $w^n_{i}=\frac{v^n_{i}}{\dt}$ are the steps of
the quasistatic evolution of the dissipated power.
If we set
\begin{equation}\label{f:gn}
w^n(x,t)=w^n_i(x)\,, \quad
g^n(x,t)= \frac{h^n_{i-1}(x)-h^n_i(x)}{\dt}\,,\qquad t\in [t^n_i, t^n_{i+1})\,
\end{equation}
and $w^n(t)\equiv w^n(\cdot,t)$, $h^n(t)\equiv h^n(\cdot,t)$, $g^n(t)\equiv g^n(\cdot,t)$,
formula (\ref{f:faradisc}) can be rewritten as
\begin{equation}\label{f:faradcon}
-\dive\left(w^n(t)\,D\gauge(Dh^n(t))\right)= g^n(t).
\end{equation}

\begin{theorem}\label{t:limitpow}
If $H_s$ is a nondecreasing function in $[0,T]$, then
the sequence
$(w^n(x,t))$ converges pointwise in $\overline{\Omega} \times [0,T]$ to
\[
w(x,t)=\frac{\partial H_s}{\partial t}(t)\,\int_{\distm(x)}^{\len^-(y)}
\chi_{\left\{h_0(y+s D\gaugem(\nor(y)))\leq H_s(t)-s\right\}}
M_{x}^-(s)\, ds\,.
\]
If $H_s$ is a nonincreasing function in $[0,T]$, then
the sequence
$(w^n(x,t))$ converges pointwise in $\overline{\Omega}\times [0,T]$ to
\[
w(x,t)=-\frac{\partial H_s}{\partial t}(t)\,\int_{\dist(x)}^{\len(y)}
\chi_{\left\{h_0(y+s D\gauge(\nor(y)))\geq H_s(t)+s\right\}}
M_{x}(s)\, ds\,.
\]
Moreover, in both cases, the sequence $(w^n)$ converges to $w$ in the strong
topology of $L^p(\Omega)$, $p\geq 1$, uniformly in $[0,T]$.
\end{theorem}

\begin{proof}
We compute the limit of the
sequence $(w^n)$ only in the case of nondecreasing external field, the other
case being similar.

We can assume, without loss of generality, that for every
$n\in\N$ the partition $P^{n+1}$ is a refinement of the
partition $P^n$. Hence, for every $t\in [0,T]$ and for every
$n\in\N$ there exists $i=i(n)\in\N$ such that $t\in [t^n_{i(n)},t^n_{i(n)+1})$.
Moreover, since $H_s\in C^1([0,T])$, we obtain
\begin{equation}\label{f:unifg}
-\frac{\partial H_s}{\partial t}(t)=\lim_{n\to \infty}
\frac{H_s(t^n_{i-1})-H_s(t^n_{i})}{\dt} \qquad
\textrm{uniformly in\ }[0,T]\,.
\end{equation}
By Proposition \ref{p:vni} we get
\[
 w^n(x,t)=
\int_{\distm(x)}^{\len^-(y)}g^n(y+s D\gaugem(\nor(y)),t)
M_{x}^-(s)\, ds\,,
\]
where $g^n$ is the function defined in (\ref{f:gn}).
Moreover we have
\begin{equation}\label{f:unifg2}
\begin{split}
h^n_{i-1} - h^n_i & =(H_s(t^n_{i-1})-H_s(t^n_{i}))
\chi_{\left\{h_0 \leq H_s(t^n_{i-1})-\distm\right\}}
\\ & +
(h_0-H_s(t^n_{i})+\distm)
\chi_{\left\{H_s(t^n_{i-1})-\distm < h_0 < H_s(t^n_{i})-\distm\right\}}\,,
\end{split}
\end{equation}
and
\[
\begin{split}
0 & \geq (h_0-H_s(t^n_{i})+\distm)
\chi_{\left\{H_s(t^n_{i-1})-\distm < h_0 < H_s(t^n_{i})-\distm\right\}}\\
& \geq \left(H_s(t^n_{i-1})-H_s(t^n_{i})\right)
\chi_{\left\{H_s(t^n_{i-1})-\distm < h_0 < H_s(t^n_{i})-\distm\right\}}\,.
\end{split}
\]
Collecting the previous information, we get
\begin{equation}\label{f:plim1}
\lim_{n\to +\infty}g^n(t) = g(t):=
-\frac{\partial H_s}{\partial t}(t)\,\chi_{\left\{h_0\leq H_s(t)-\distm\right\}}\,,
\end{equation}
and $\lim_{n\to +\infty}\|g^n(t)-g(t)\|_{L^p(\Omega)}=0$ uniformly in $[0,T]$, for
every $p\geq 1$. In addition, from (\ref{f:unifg}) and (\ref{f:unifg2}),
for every regular point $x$ of $\distm$ and setting $\{y\}=\Pi(x)$, we get
\begin{equation}\label{f:plim3}
\lim_{n\to +\infty}g^n(y+sD\gaugem(\nor(y)), t)= g(y+sD\gaugem(\nor(y)), t)\qquad
\textrm{a.e.\ }s\in [\distm(x), \len^-(y)]\,.
\end{equation}
Finally,
by (\ref{f:plim3}), (\ref{f:estiM}) and the Dominated Convergence Theorem, we conclude that
\begin{equation}\label{f:plim2}
\begin{split}
&\lim_{n \to +\infty} w^n(x,t)=\\
& -\frac{\partial H_s}{\partial t}(t)\,\int_{\distm(x)}^{\len(y)}
\chi_{\left\{h_0(y+s D\gaugem(\nor(y)))\leq H_s(t)-s\right\}}
M_{x}^-(t)\, ds\,:=w(x,t)\,.
\end{split}
\end{equation}
The last part of the theorem follows from the fact that
\begin{equation}\label{f:plimw}
{\|w^n(t)-w(t)\|}_{L^p(\Omega)}\leq
C_0\, {\|g^n(t)-g(t)\|}_{L^p(\Omega)}
\end{equation}
for some positive constant $C_0$.
Namely, for a given $t\in [0,T]$ let us define
\[
\varphi^n(x, s) := g^n(y+sD\gaugem(\nor(y)), t) - g(y+sD\gaugem(\nor(y)), t),
\]
for $x\in\Omega\setminus\overline{\Sigma}$,
$\proj(x)=\{y\}$, and $s\in [0,\len^-(y)]$.
We remark that, if $x\in\Omega\setminus\overline{\Sigma}$ and
$y$ is the projection of $x$,
then $\varphi^n(x,\cdot) = \varphi^n(y,\cdot)$.

By the very definition of $w^n$, $w$ and $M_x$, H\"older's inequality and
(\ref{f:estiM}) we have that
\[
{\|w^n(t)-w(t)\|}^p_{L^p(\Omega)}\leq
C \int_{\Omega} \left(
\int_{\distm(x)}^{\len^-(x)} |\varphi^n(x,s)|^p
\prod_{i=1}^{n-1}\frac{1-s\curvgm_i(x)}{1-\distm(x)\curvgm_i(x)}
\, ds
\right)dx =: I
\]
where $\len^-(x):= \len^-(y)$ when $\proj^-(x) = \{y\}$ and
$C := M_0^{p-1}\, (\max\distm)^{1/p'}$.
Using the change of variables theorem (see \cite[Thm.~7.1]{CMf}) we
have that
\[
\begin{split}
I = {} &
C \int_{\partial\Omega} \gaugem(\nor(y)) \left[\int_0^{\len^-(y)}
\left(\int_{\sigma}^{\len^-(y)} |\varphi^n(y,s)|^p
\prod_{i=1}^{n-1}\frac{1-s\curvgm_i(y)}{1-\sigma\curvgm_i(y)}
\, ds \right) \right.\cdot
\\ & {} \cdot
\left.\prod_{i=1}^{n-1} (1-\sigma\, \curvg^-_i(y))\, d\sigma\right]
d\haus(y)
\\ = {} &
C \int_{\partial\Omega} \gaugem(\nor(y)) \left[\int_0^{\len^-(y)}
\left(\int_{\sigma}^{\len^-(y)} |\varphi^n(y,s)|^p
\prod_{i=1}^{n-1}(1-s\curvgm_i(y))
\, ds \right) \, d\sigma\right]
d\haus(y)
\\ \leq {} &
C' \int_{\partial\Omega} \gaugem(\nor(y)) \left(\int_0^{\len^-(y)}
|\varphi^n(y,s)|^p
\prod_{i=1}^{n-1}(1-s\curvgm_i(y))
\, ds \right)
d\haus(y),
\end{split}
\]
where $C' := C\, \max\distm$.
Using again the change of variables theorem we finally get
\[
I \leq C' \int_{\Omega} |g^n(x,t)-g(x,t)|^p\, dx,
\]
hence (\ref{f:plimw}) follows.
\end{proof}

\begin{rk}
Notice that if $\distm(x)>\lambda^n_i(y)$, then $v^n_i(x)=0$, while
$h^n_i(x)=c^n_i-\distm(x)$ if $\distm(x)\leq \lambda^n_i(y)$. Hence we have that
$w^n(x,t) D\gauge(Dh^n(x,t))=w^n(x,t) D\gauge(-D\distm(x))$ a.e.\ in $\Omega$ and
for every $t\in [0,T]$.
A passage to the limit in (\ref{f:faradcon})
leads to
\begin{equation}\label{f:faradfin}
-\dive(w(x,t)D\gauge(Dh(x,t)))=- \frac{\partial H_s}{\partial t}
\,\chi_{\left\{h_0\leq H_s(t)-\distm\right\}}\,,
\end{equation}
in the sense of distributions in $\Omega$ for every $t\in [0,T]$.
Recalling that, by construction, $\frac{\partial H_s}{\partial t}
\,\chi_{\left\{h_0\leq H_s(t)-\distm\right\}}= \frac{\partial h}{\partial t}(x,t)$, and
comparing (\ref{f:faradfin}) with the Faraday's law, we conclude
that the electric field induced by $h$ inside the superconductor is
$\vec{E}(x,t)=w(x,t)D\gauge(Dh(x,t))$.
\end{rk}

\section{Further remarks}

As a consequence of the Theorem \ref{t:limitsol} we obtain some
detailed information about the macroscopic behavior of the
penetrating magnetic field.

\begin{center}
\begin{figure}
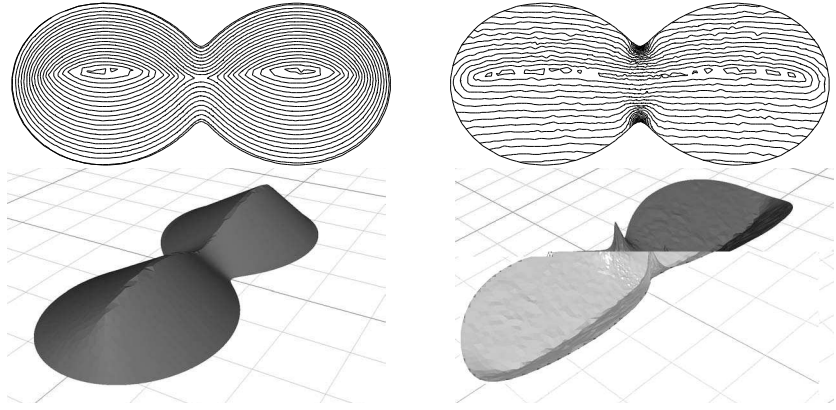

\includegraphics[width=5cm]{\picturename{ECdistbw-crop}}\qquad
\includegraphics[width=5cm]{\picturename{ECvbw-crop}}
\includegraphics[width=5cm]{\picturename{ECdist3-crop}} \qquad
\includegraphics[width=5cm]{\picturename{ECv3-crop}}
\caption{Level sets and 3D plot of the full penetrated magnetic
field  and of the dissipation for $\Omega$ Cassini's Egg, $K$ ellipse not centered at the origin,
$h_0=0$}\label{BB3}
\end{figure}
\end{center}

\noindent\textit{Saturation time (Full penetration time).}
The experiments show that, if the superconductor is posed in an
increasing external field $H_s$, say $H_s(t) \to +\infty$ as
$t \to +\infty$, after a while the internal magnetic field $h(x,t)$ differs
from $H_s(t)$ for a stationary amount $u(x)$.
In our model we obtain that if $\tau \geq 0$ is such that
$H_s(\tau)=\sup_\Omega (h_0+\distm)$, then for $t\geq \tau$
$h(x,t)-H_s(t)=-\distm(x)$.

\begin{center}
\begin{figure}
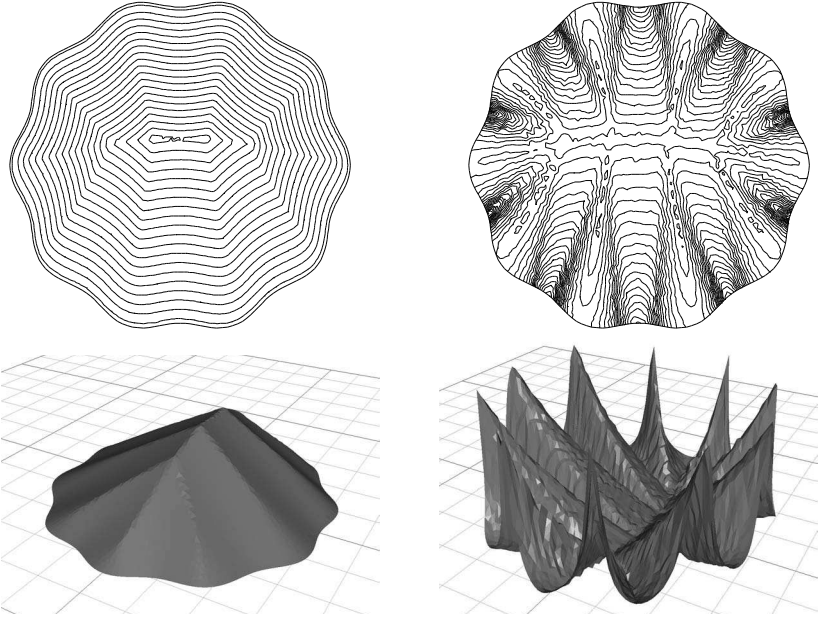

\includegraphics[width=4.5cm]{\picturename{EFdistbw-crop}}\qquad\qquad
\includegraphics[width=4.5cm]{\picturename{EFvbw-crop}}
\includegraphics[width=5cm]{\picturename{EFdist3-crop}}\qquad
\includegraphics[width=5cm]{\picturename{EFv3-crop}}
\caption{Level sets and 3D plot of the full penetrated magnetic
field and of the dissipation for $\Omega$ perturbed disk, $K$ ellipse not centered at the origin,
$h_0=0$.}
\label{BB4}
\end{figure}
\end{center}

In the simpler case when $H_s(t)=-at$, $a>0$, and
$h_0=0$ we obtain that at time $\tau=\frac{\max_\Omega \dist}{a}$
the magnetic field is penetrated in the whole superconductor and
$h(x,t)=\dist(x)+H_s(t)$ for $t\geq \tau$.
Figures \ref{BB3} and  \ref{BB4} picture
the plots of the magnetic field and of the dissipation in the
stationary configuration for anisotropic materials with a
complicated geometry. We underline that
the numerical computation based on the representation formula
(\ref{f:plim2}) shows the presence of an hight dissipation near
the points of the boundary with negative curvature, according to the
experimental observations.

\noindent\textit{Boundary of penetration profile.}
If the evolution starts with $h_0=0$ (no initial magnetic
field inside $\Omega$) and the external field is nondecreasing,
then at time $t$ the induced magnetic
field is penetrated in $\Omega$ only in the region
$\{x\in \Omega \colon\ \distm(x)\leq H_s(t)\}$. This gives an
exact macroscopic description of the boundary of penetration profile,
which is given by the level sets of the Minkowski distance from
the boundary of $\Omega$.

\begin{center}
\begin{figure}
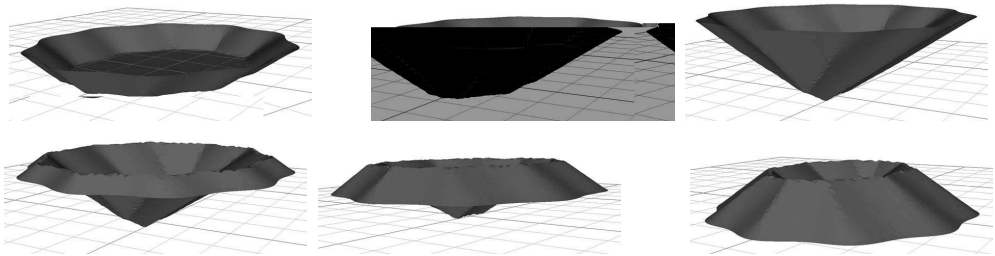

\includegraphics[width=4cm]{\picturename{hyst0-crop}}\qquad
\includegraphics[width=4cm]{\picturename{hyst1-crop}}
\includegraphics[width=4cm]{\picturename{hyst2-crop}}\qquad
\includegraphics[width=4cm]{\picturename{hyst3-crop}}
\includegraphics[width=4cm]{\picturename{hyst4-crop}} \qquad
\includegraphics[width=4cm]{\picturename{hyst5-crop}}
\caption{$\Omega$ perturbed disk, $K$ ellipse centered at the origin. Profile of the internal magnetic field,
starting from $h_0=0$, during a loop of the external field.}\label{BB5}
\end{figure}
\end{center}

\noindent\textit{Hysteresis phenomenon.}
Assume that the evolution starts with $h_0=0$ (no initial magnetic
field inside $\Omega$) and that the external field $H_s(t)$ is increasing
in [0,T], then at time $T$ the induced magnetic
field is $h(x,T)=(-\distm(x)+H_s(t))_+$.
Starting from this configuration,
in [T,2T] the superconductor is subject to the external field
$H_s(2T-t)$. Then we have
$h(x,2T)=\min\{(-\distm(x)+H_s(T))_+, \dist(x)+H_s(0)\}$ which in general
is not zero. This corresponds with the experimental observation of
an hysteresis phenomenon for these hard superconductors  after an external field loop.
\begin{center}
\begin{figure}
\includegraphics[width=4cm]{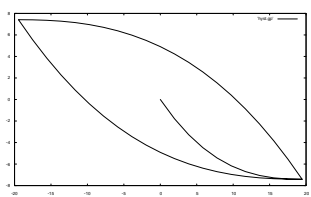}
\caption{Magnetic hysteresis loop corresponding to the previous process.}
\end{figure}
\end{center}

\noindent\textit{On the validity of the Badi\'a--L\'opez formalism for anisotropic
materials.} In \cite{Fisher} Fisher doubts the validity of the variational
formulation proposed by Bad\'ia and Lopez in the case of anisotropic
materials. He declares that \textit{``one cannot be sure that the suggested
least action principle chooses the Current--Voltage Characteristics correctly.
Of course, if this principle is proved, it would be a very important tool
to determine the Current--Voltage Characteristics via the known restriction
region''}. Our results state that, if the power--law approximation for the dissipation
is assumed to be reliable, then the Current--Voltage Characteristics
induced by this approximation is the one selected by the minimization procedure.

%

\begin{thebibliography}{10}

\bibitem{BLd}
{A.}~Bad\'\i a and {C.} L\'opez, \emph{Vector magnetic hysteresis of hard
  superconductors}, Phys.\ Rev.~B \textbf{65} (2001), 104514.

\bibitem{BLc}
{A.}~Bad\'\i a and {C.} L\'opez, \emph{Horizons in Superconductivity Research},
  ch.~Minimal model for the topology of the critical state in hard
  superconductors, Nova Science Publishers, 2003.

\bibitem{BLb}
{A.}~Bad\'\i a and {C.} L\'opez, \emph{Electric field in hard superconductors
  with arbitrary cross section and general critical current law}, J.\ Appl.\
  Phys. \textbf{95} (2004), 8035--8040.

\bibitem{BaCD}
{M.} Bardi and {I.} Capuzzo~Dolcetta, Optimal Control and Viscosity Solutions
  of {H}amilton-{J}acobi-{B}ellman Equations, Systems \& Control: Foundations
  \& Applications, Birkh\"auser, Boston, 1997.

\bibitem{BaPrb}
{J.W.} Barrett and {L.} Prigozhin, \emph{Sandpiles and superconductors: dual
  variational formulations for critical-state problems}, preprint, 2005.

\bibitem{BaPr}
{J.W.} Barrett and {L.} Prigozhin, \emph{{B}ean's critical-state model as the
  $p\to\infty$ limit of an evolutionary $p$-{L}aplacian equation}, Nonlinear
  Anal. \textbf{42} (2000), 977--993.

\bibitem{Bean}
{C.P.} Bean, \emph{Magnetization of hard superconductors}, Phys.\ Rev.\ Letters
  \textbf{8} (1962), 250--253.

\bibitem{BKR}
{K.V.} Bhagwat, {D.} Karmakar, and {G.} Ravikumar, \emph{Critical state model
  with anisotropic critical current density}, J.\ Phys.: Condens.\ Matter
  \textbf{15} (2003), 1325--1337.

\bibitem{Brai}
{A.} Braides, $\Gamma$-convergence for beginners, Oxford University Press, New
  York, 2002.

\bibitem{Bran}
{E.H.} Brandt, \emph{Electric field in superconductors with rectangular cross
  section}, Phys.\ Rev.~B \textbf{52} (1995), 15442--15457.

\bibitem{Butt}
{G.} Buttazzo, Semicontinuity, relaxation and integral representation in the
  calculus of variations, Pitman Res.\ Notes Math.\ Ser., vol. 207, Longman
  Scientific and Technical, Harlow, U.K., 1989.

\bibitem{CaSi}
{P.} Cannarsa and {C.} Sinestrari, Semiconcave functions, {H}amilton-{J}acobi
  equations and optimal control, Progress in Nonlinear Differential Equations
  and their Applications, vol.~58, Birkh\"auser, Boston, 2004.

\bibitem{Cha}
{S.J.} Chapman, \emph{A Hierarchy of Models for Type-{II} Superconductors},
  SIAM Rev. \textbf{42} (2000), 555--598.

\bibitem{CMf}
{G.} Crasta and {A.} Malusa, \emph{The distance function from the boundary in a
  {M}inkowski space}, to appear in Trans.\ Amer.\ Math.\ Soc.,
  oai:ar{X}iv:math.AP/0612226.

\bibitem{CMg}
{G.} Crasta and {A.} Malusa, \emph{On a system of partial differential
  equations of {M}onge-{K}antorovich type}, to appear in J.\ Differential
  Equations, oai:ar{X}iv:math.AP/0612227.

\bibitem{DM}
{G.} Dal~Maso, An introduction to $\Gamma$--convergence, Birkh\"auser, Boston,
  1993.

\bibitem{Fisher}
{L.M.} Fisher and {V.A.} Yampol'skii, \emph{Comment on ``Critical statetheory
  for non-parallel flux line lattices in type-{II} superconductors''},
  oai:ar{X}iv:cond-mat/0201286, 2002.

\bibitem{GNP}
{A.} Garroni, {V.} Nesi, and {M.} Ponsiglione, \emph{Dielectric breakdown:
  optimal bounds}, Proc.\ R.\ Soc.\ Lond.\ A \textbf{457} (2001), 2317--2335.

\bibitem{Gold}
{S.I.} Goldberg, Curvature and homology, Academic Press, New York, 1970.

\bibitem{LN}
{Y.Y.} Li and {L.} Nirenberg, \emph{The distance function to the boundary,
  {F}insler geometry and the singular set of viscosity solutions of some
  {H}amilton--{J}acobi equations}, Commun.\ Pure Appl.\ Math. \textbf{58}
  (2005), 85--146.

\bibitem{Li}
{P.L.} Lions, Generalized solutions of {H}amilton-{J}acobi equations, Pitman,
  Boston, 1982.

\bibitem{Rock}
{R.T.} Rockafellar, Convex Analysis, Princeton Univ.\ Press, Princeton, NJ,
  1970.

\bibitem{Sch}
{R.} Schneider, Convex bodies: the {B}runn--{M}inkowski theory, Cambridge
  Univ.~Press, Cambridge, 1993.

\bibitem{Wa}
{W.} von Wahl, \emph{Estimating $\nabla u$ by div${}\,u$ and curl${}\,u$},
  Math.\ Methods Appl.\ Sci. \textbf{15} (1992), 123--143.

\bibitem{YLZ}
{H.-M.} Yin, {B.Q.} Li, and {J.} Zou, \emph{A degenerate evolution system
  modeling {B}ean's critical-state type-{II} superconductors}, Discrete
  Contin.\ Dynam.\ Systems \textbf{8} (2002), 781--794.

\end{thebibliography}

\end{document}